\documentclass[pdflatex,sn-mathphys-num]{sn-jnl}

\usepackage{amsmath}
\usepackage{amssymb}
\usepackage{amsfonts}
\usepackage{amsthm}
\usepackage{comment}
\usepackage{pgfplots}
\usepackage{graphicx}
\usepackage{amsmath}
\usepackage[T1]{fontenc}
\usepackage{mathrsfs}
\usepackage{latexsym}
\usepackage{amssymb}
\usepackage{graphicx}
\usepackage{float}
\usepackage{extarrows}
\usepackage{epstopdf}
\usepackage{caption}
\usepackage{subcaption}
\usepackage{amsmath}
\usepackage{wrapfig}
\usepackage{lipsum}
\usepackage{array}
\graphicspath{ {./Figures/}}

\newtheorem{remark}{Remark}
\newtheorem{theorem}{Theorem}

\newtheorem{lemma}{Lemma}
\newtheorem{corollary}{Corollary}

\begin{document}

\title[]{A novel efficient structure-preserving exponential integrator for Hamiltonian systems}


\author[1]{\fnm{Pan} \sur{Zhang}}\email{zhangp273@mail2.sysu.edu.cn}

\author[2]{\fnm{Fengyang} \sur{Xiao}}\email{xiaofy5@mail2.sysu.edu.cn}
\author*[1]{\fnm{Lu} \sur{Li}}\email{lilu86@mail.sysu.edu.cn}


\affil[1]{\orgname{Sun Yat-sen University}, \orgaddress{\city{Zhuhai}, \postcode{519082}, \country{China}}}

\affil[2]{\orgname{Duke University}, \orgaddress{\city{North Carolina}, \postcode{27705}, \country{United States}}}

\abstract{We develop an efficient linearly implicit time-reversible exponential integrator for semilinear Hamiltonian systems with polynomial nonlinearities. The method combines the variation of constants formula with  Kahan-type discretization of the nonlinear term. For quadratic vector fields,  we prove second-order convergence in finite dimensions and derive an exact stepwise energy-increment identity implying a precise notion of near  conservation of energy. For higher-degree polynomial Hamiltonians, we employ polarization to obtain a linearly implicit multistep exponential extension and discuss its well-posedness, stability and convergence; we further provide a stepwise energy characterization that quantifies its energy drift in terms of the numerical increments. Numerical experiments on the H\'enon-Heiles system, a continuum Fermi-Pasta-Ulam model, and the two-dimensional Zakharov-Kuznetsov equation corroborate the theoretical accuracy and demonstrate competitive accuracy-cost tradeoffs compared with representative symmetric energy-preserving exponential integrators.}

\keywords{Hamiltonian system,  linearly implicit, exponential integrato, symmetric, near conservation of energy}



\maketitle

\section{Introduction}
This study is concerned with semilinear systems of the form:
\begin{equation}\label{semilinear ODE}
	\dot{x}(t) = Ax(t) + f(x(t)), \quad x(t_0) = x_0,
\end{equation}
where the matrix $A$ often represents a spatial discretization of a linear unbounded differential operator, and its eigenvalues have large negative real parts or are purely imaginary with substantial magnitude. The nonlinear term $f$ satisfies the Lipschitz condition. Such systems appear in a wide range of applications, including mechanics \cite{ambati2015review}, quantum physics \cite{kato1987nonlinear}, charged-particle dynamic \cite{Cary} and discretizations of partial differential equations (PDEs) \cite{Celledoni2008}. For these types of systems, exponential integrators are commonly employed due to their ability to efficiently handle the linear part through matrix exponentials, allowing for larger step sizes and improved accuracy in stiff and oscillatory problems \cite{MR2652783}.

Semilinear systems \eqref{semilinear ODE} often exhibit  Hamiltonian structures like:
\begin{equation}\label{semilinear HODE}
	\begin{split}
		\dot{x}(t) &= Q \nabla H(x(t)),\\ 
		x(t_0) &= x_0, \quad H(x(t))=\frac{1}{2}x^\mathsf{T} Mx+U(x),
	\end{split}
\end{equation}
i.e., $A=QM$ and $f(x)=Q\nabla U(x)$, where $Q$ is a skew-symmetric matrix, $H$ is the Hamiltonian energy function containing a quadratic part  and a higher-degree part \cite{Celledoni2012}. For such systems, the Hamiltonian remains invariant along the exact solution trajectory. Moreover, when $Q$ is a canonical skew-symmetric matrix, the system preserves the symplectic two-form. Maintaining these geometric structures is essential for the long-time numerical integration of such systems, as it ensures accuracy and stability \cite{hairer2006geometric}.  To address these  needs, various structure-preserving integrators have been proposed, including symmetric exponential integrators \cite{Celledoni2008,li2022new}, symplecticity-related exponential integrators  \cite{wu2012explicit,Wuxinyuan2018,mei2019exponential,jiang2022high,sato2024high}, volume-preserving exponential integrators \cite{wang2019volume} and energy-preserving exponential integrators \cite{Wuxinyuan2016,fu2022high}.  Extensive studies have demonstrated that these specialized methods deliver significantly enhanced long-term performance compared to conventional high-order schemes.

In recent years, energy-preserving exponential integrators have been widely studied due to their effectiveness in accurately capturing the energy dynamics of Hamiltonian systems \cite{Wuxinyuan2016,MEI2022110822,jiang2020linearly,MEI2021110429,DENG2022166,Ju2022,mei2022energy,li2022energy}. A prominent example is the exponential averaged vector field (EAVF) method \cite{Wuxinyuan2016} which combines the discrete gradient approach with exponential integrator. This hybrid strategy has inspired the development of a broad class of exponential energy-preserving schemes. Solving such fully implicit methods requires tackling nonlinear systems at each time step, which can be  computationally demanding. As an alternative, linearly implicit structure-preserving methods  have been explored. These methods offer similar stability to fully implicit schemes while requiring only the solution of a single linear system per time step, making them more computationally efficient.
One class of such methods is based on the Scalar Auxiliary Variable (SAV) framework \cite{shen2019new,liu2020,jiang2020linearly}, which conserves a modified energy for semilinear Hamiltonian systems but lacks symmetry. Symmetry is also an important geometric property, as it enforces discrete time-reversibility, yields an even-power error expansion, and often improves long-time behavior \cite{hairer2006geometric}. 
A number of further developments along these lines can be found in \cite{jiang2022high,sato2024high,Ju2022,Xu2022,fu2022arbitrary}.   Another category involves  multiple-point methods that portion out  the
nonlinearity \cite{matsuo2001dissipative,li2022new}. These methods have the advantage of preserving both a modified energy and symmetry, but they  may suffer from stability concerns when compared to Kahan methods \cite{eidnes2019linearly}, and a systematic convergence theory for linearly implicit time-reversible exponential integrators, especially in the higher-degree polynomial setting, is still limited \cite{li2022new}. {\color{black}
	These observations motivate the design of stable time-reversible and linearly implicit exponential integrators, together with a rigorous convergence analysis. We therefore incorporate the Kahan discretization for the nonlinear term, since it is linearly implicit, time-reversible and  known to preserve a modified energy and a modified measure \cite{MR3452144,celledoni2012geometric}.}

{\color{black}To treat the stiff linear dynamics in \eqref{semilinear ODE}, we embed Kahan update into the variation of constants exponential framework, leading to symmetric schemes that require only a linear solve per step and admit a quantitative near  conservation of energy. This coupling is, however, nontrivial: the exponential update changes the algebraic structure of the original Kahan discretization, so its classical modified-energy preservation mechanism does not carry over in general, and the energy conservation analysis must be rebuilt.
	Our main contributions in this work are as follows.
	First, for quadratic vector fields (cubic Hamiltonians), we obtain a one-step linearly implicit symmetric exponential scheme and prove second-order convergence on finite time intervals.
	Second, for higher-degree polynomial Hamiltonians, we derive a linearly implicit time-reversible multistep exponential extension via polarization and establish  solvability and stability analysis.
	Third, we derive explicit stepwise energy increment identities, which yield computable energy-drift estimates in terms of numerical increments and  explain the favorable long-time behavior observed in our experiments.
}

The remainder of this paper is organized as follows. In Section \ref{Method}, {\color{black}we introduce the proposed linearly implicit  exponential integrator for cubic Hamiltonians,  discuss its key properties, including symmetry, convergence, near  preservation of energy, and its generalization to systems with higher-degree polynomial Hamiltonians.} In Section \ref{Experiment}, we present numerical experiments to demonstrate the stability, computational efficiency, and long-term behavior of the proposed method using different types of differential equations. Finally, conclusions and potential directions for future research are discussed in the last section.

\section{The new linearly implicit  structure-preserving symmetric method}\label{Method}

Before proposing the new exponential integrator, we will introduce Kahan's method and exponential integrators since our new scheme is created based on them. For system \eqref{semilinear HODE} with a  cubic Hamiltonian,   Kahan's method has the form 
\begin{equation*}
	x_{n+1}= Q(-\frac{1}{2}\nabla H(x_{n})+2\nabla H(\frac{x_{n}+x_{n+1}}{2})-\frac{1}{2}\nabla H(x_{n+1})),
\end{equation*}
and it is shown to preserve a modified energy as well as a measure in \cite{celledoni2012geometric}. In particular, if the energy function is homogeneous, the modified energy can be written as $\bar{H}=C(x_n,x_n,x_{n+1})$, where $C$ is a symmetric trilinear form defined by the Hessian of the energy function $H$
\begin{equation*}
	C(x,y,z) = \frac{1}{6}x^\mathsf{T} H^{\prime\prime}(y)z.
\end{equation*}
When $H$ is nonhomogeneous, we can always add an extra auxiliary variable to transfer the problem to a homogeneous setting.
In addition, Kahan's method is linearly implicit. These properties guarantee the superior numerical behavior  of Kahan's method during long-time integration \cite{eidnes2019linearly}.

For a semilinear system \eqref{semilinear ODE}, we consider the following variation of constants formula on interval $[t_{0},t_{0}+h]$  
\begin{equation*}
	x(t_{0}+h)=e^{hA} x(t_{0})+h\int_{0}^{1} e^{(1-\tau)hA}f(x(t_{0}+\tau h))\mathrm{d} \tau.
\end{equation*}
The key idea of constructing exponential integrator is to integrate the linear part exactly meanwhile finding an approximation to the integration of the nonlinear part. For example, the exponential Euler method is obtained by approximating the nonlinear term by $f(x(t_0))$ \cite{MR2652783}, thereby having the form
\begin{equation*}\label{Eeuler}
	x_{n+1}=e^{hA} x_{n} +h \phi(hA)f(x_{n}),
\end{equation*}
where $\phi(z):= \frac{e^{z}-1}{z}$;  the exponential averaged vector field (EAVF) method is obtained by approximating the nonlinear term by the discrete gradient \cite{Wuxinyuan2016}, possessing the form
\begin{equation*}\label{EAVF}
	x_{n+1}=e^{hA} x_{n} +h \phi(hA)\int_0^1f(\xi x_{n}+(1-\xi)x_{n+1})d\xi.
\end{equation*}

\subsection{{\color{black}The one-step EKahan scheme for cubic Hamiltonians}}
For system \eqref{semilinear HODE} with a cubic Hamiltonian, we propose a linearly implicit exponential integrator by approximating the nonlinear term using the above mentioned Kahan discretization, denoted by EKahan 
\begin{equation}\label{EKahan scheme}
	x_{n+1}=e^{hA} x_{n} +h \phi(hA)(-\frac{1}{2}f(x_{n})+2f(\frac{x_{n}+x_{n+1}}{2})-\frac{1}{2}f(x_{n+1})).
\end{equation}
Before presenting the properties of scheme \eqref{EKahan scheme}, we give the following notations and assumptions:
\paragraph{Norms and growth constants}
Let $\|\cdot\|$ be an induced matrix norm and the associated vector norm.
For a fixed final time $T>0$, define
\begin{equation*}
	M_E=\max_{0\le t\le T}\|e^{tA}\|,\qquad
	M_\phi=\max_{0\le t\le T}\|\phi(tA)\|.
\end{equation*}
\paragraph{Lipschitz assumptions}
For system \eqref{semilinear HODE}, we assume that there exists a convex bounded set $\mathcal D\subset\mathbb R^d$ such that both
the exact solution $\{x(t):t\in[0,T]\}$ and the numerical solution $\{x_n:0\le nh\le T\}$
remain in $\mathcal D$, and that $f$ is Lipschitz on $\mathcal D$ with constant $L$.
\begin{lemma}\label{EKahan order Lemma}
	{\color{black}Assuming that $f$ is a quadratic vector field and that the exact solution of \eqref{semilinear ODE} is sufficiently smooth on $[t_n,t_{n+1}]$, then the EKahan scheme \eqref{EKahan scheme} has local truncation error $O(h^3)$.}
\end{lemma}
\begin{proof}
	{\color{black}For quadratic $f$, we have the exact Taylor formulas
			\begin{align*}
				f(x_{n+1}) & = f(x_n) + f'(x_n)\Delta_x + \tfrac12 f''(x_{n})(\Delta_x,\Delta_x),\\
				f\left(\tfrac{x_n+x_{n+1}}{2}\right) & =f(x_n)+\tfrac12 f'(x_n)\Delta_x+\tfrac18 f''(x_{n})(\Delta_x,\Delta_x),
			\end{align*}
		where $\Delta_x=x_{n+1}-x_n$ and $f''$ is constant. Hence, we have
		\begin{equation}\label{eq:quad-identity}
			-\tfrac12 f(x_n)+2f\!\left(\tfrac{x_n+x_{n+1}}{2}\right)-\tfrac12 f(x_{n+1})
			= f(x_n)+\tfrac12 f'(x_n)\Delta_x.
		\end{equation}
		Using \eqref{eq:quad-identity}, the EKahan scheme \eqref{EKahan scheme} can be rewritten as
		\begin{equation}\label{eq:ekahan-rewrite}
			x_{n+1}=e^{hA}x_n + h\phi(hA)\Bigl(f(x_n)+\tfrac12 f'(x_n)\Delta_x\Bigr).
		\end{equation}
		Since $e^{hA}=I+h\phi(hA)A$ and $A$ commutes with $\phi(hA)$, \eqref{eq:ekahan-rewrite} becomes
		\[
		x_{n+1}=x_n + h\phi(hA)g(x_n) + \tfrac{h}{2}\phi(hA)f'(x_n)\Delta_x,
		\quad \text{with}\quad g(x)=Ax+f(x),
		\]
		and therefore we have 
		\begin{equation}\label{eq:Delta-eq}
			\Bigl(I-\tfrac{h}{2}\phi(hA)f'(x_n)\Bigr)\Delta_x = h\phi(hA)g(x_n).
		\end{equation}
		For sufficiently small $h$, the matrix in \eqref{eq:Delta-eq} is invertible, and we may write
		\begin{equation}\label{eq:Delta-sol}
			\Delta_x = h\Bigl(I-\tfrac{h}{2}\phi(hA)f'(x_n)\Bigr)^{-1}\phi(hA)g(x_n).
		\end{equation}
		Using the series expansion $\phi(hA)=I+\tfrac{h}{2}A+O(h^2)$ and the Neumann expansion
		\[
		\Bigl(I-\tfrac{h}{2}\phi(hA)f'(x_n)\Bigr)^{-1}
		= I+\tfrac{h}{2}\phi(hA)f'(x_n)+O(h^2),
		\]
		we obtain from \eqref{eq:Delta-sol} that 
		\begin{align*}
			\Delta_x
			= h\,g(x_n) + \tfrac{h^2}{2}\bigl(A+f'(x_n)\bigr)g(x_n) + O(h^3).
		\end{align*}
		Consequently, we have 
		\begin{equation}\label{eq:method-expansion}
			x_{n+1}=x_n + h\,g(x_n) + \tfrac{h^2}{2}\bigl(A+f'(x_n)\bigr)g(x_n) + O(h^3).
		\end{equation}
		On the other hand, the exact solution satisfies
		\begin{equation}\label{eq:exact-expansion}
			x(t_{n+1})=x(t_n)+h\,g(x(t_n))+\tfrac{h^2}{2}g'(x(t_n))g(x(t_n))+O(h^3),
		\end{equation}
		where $g'(x)=A+f'(x)$. Taking $x_n=x(t_n)$ and comparing \eqref{eq:method-expansion} with \eqref{eq:exact-expansion} yields $x_{n+1}-x(t_{n+1}) = O(h^3)$.
	} 
\end{proof}

\begin{theorem}\label{thm:global-error}
	{\color{black}Assume $f$ satisfies the standard Lipschitz condition and the local truncation error satisfies $\|\tau_{n+1}\|\le C_{\rm loc}h^3$ for all $0\le nh\le T$, which is true by Lemma~\ref{EKahan order Lemma}.
		Then there exist constants $C_1,C_2>0$, depending only on $\|A\|$, $L$, $T$, $M_E$, $M_\phi$ and $C_{\rm loc}$, such that
		\[
		\|x(t_n)-x_n\|\le e^{C_1T}\Bigl(\|x(t_0)-x_0\|+C_2h^2\Bigr),\quad 0\le nh\le T,
		\]
		if the step size satisfies
		\begin{equation}\label{eq:h0}
			h\le h_0:=\frac{1}{3LM_\phi}.
		\end{equation}
		In particular, if $x_0=x(t_0)$, then $\max\limits_{0\le nh\le T}\|x(t_n)-x_n\|=O(h^2)$}.
\end{theorem}

\begin{proof}
	Set $\varepsilon_n=x(t_n)-x_n$ and 
	$\Psi(u,v)=-\tfrac12 f(u)+2f\!\Bigl(\tfrac{u+v}{2}\Bigr)-\tfrac12 f(v)$.
	Inserting the exact solution into \eqref{EKahan scheme} yields
	\begin{equation}\label{eq:exact-scheme2}
		x(t_{n+1})=e^{hA}x(t_n)+h\phi(hA)\Psi(x(t_n),x(t_{n+1}))+\tau_{n+1},
		\quad \|\tau_{n+1}\|\le C_{\rm loc}h^3.
	\end{equation}
	Subtracting \eqref{EKahan scheme} from \eqref{eq:exact-scheme2} gives the error recursion
	\begin{equation}\label{eq:error-rec2}
		\varepsilon_{n+1}=e^{hA}\varepsilon_n+h\phi(hA)\delta_{n+1,n}+\tau_{n+1},
	\end{equation}
	where $\delta_{n+1,n}=\Psi(x(t_n),x(t_{n+1}))-\Psi(x_n,x_{n+1})$.
	Since $\mathcal{D}$ is convex and $f$ is Lipschitz on $\mathcal{D}$, we have 
	\begin{align*}
		\|\delta_{n+1,n}\|
		&\le \tfrac12\|f(x(t_n))-f(x_n)\|
		+2\Bigl\|f\!\Bigl(\tfrac{x(t_n)+x(t_{n+1})}{2}\Bigr)-f\!\Bigl(\tfrac{x_n+x_{n+1}}{2}\Bigr)\Bigr\|\\
		&+\tfrac12\|f(x(t_{n+1}))-f(x_{n+1})\| \\
		&\le \tfrac32 L(\|\varepsilon_n\|+\|\varepsilon_{n+1}\|).
	\end{align*}
	Taking norms in \eqref{eq:error-rec2} and using $\|\phi(hA)\|\le M_\phi$ yields
	\begin{equation}\label{eq:ineq-alpha}
		\|\varepsilon_{n+1}\|
		\le \|e^{hA}\|\,\|\varepsilon_n\|
		+\tfrac32 LhM_\phi\bigl(\|\varepsilon_n\|+\|\varepsilon_{n+1}\|\bigr)
		+C_{\rm loc}h^3.
	\end{equation}
	Denoting by $\alpha=\tfrac32 LhM_\phi$, moving the term $\alpha\|\varepsilon_{n+1}\|$ to the left in \eqref{eq:ineq-alpha} gives
	\begin{equation*}
		(1-\alpha)\|\varepsilon_{n+1}\|
		\le (\|e^{hA}\|+\alpha)\|\varepsilon_n\| + C_{\rm loc}h^3,
	\end{equation*}
	and hence we have 
	\begin{equation}\label{eq:stab-basic}
		\|\varepsilon_{n+1}\|
		\le q(h)\,\|\varepsilon_n\| + \frac{C_{\rm loc}}{1-\alpha}h^3,
		\quad \text{with} \quad 
		q(h)=\frac{\|e^{hA}\|+\alpha}{1-\alpha}.
	\end{equation}
	Using $e^{hA}=I+\int_0^h Ae^{sA} \mathrm{d}s$ and $\|e^{sA}\|\le M_E$ for $0\le s\le h$, we obtain
	\begin{equation*}
		\|e^{hA}\|\le 1+h\|A\|M_E.
	\end{equation*}
	Moreover, under \eqref{eq:h0} we have $\alpha\le \tfrac12$, so that $(1-\alpha)^{-1}\le 1+2\alpha$ and $2\alpha^2\le \alpha$.
	Therefore, we get
	\begin{align*}
		q(h)
		=\frac{\|e^{hA}\|+\alpha}{1-\alpha}
		\le (\|e^{hA}\|+\alpha)(1+2\alpha) \le 1+Ch
	\end{align*}
	with $C=\|A\|M_E + 3(1+M_E)L M_\phi$.
	In addition, since $\alpha\le\tfrac12$ we have $(1-\alpha)^{-1}\le 2$, and thus
	\[
	\frac{C_{\rm loc}}{1-\alpha}h^3 \le 2C_{\rm loc}h^3.
	\]
	Substituting these bounds into \eqref{eq:stab-basic} gives
	\begin{equation}\label{eq:stab-MI+Ch}
		\|\varepsilon_{n+1}\|\le (1+Ch)\|\varepsilon_n\|+2C_{\rm loc}h^3,
		\qquad 0\le nh\le T,\ \ h\le h_0.
	\end{equation}
	Applying the discrete Gronwall inequality to \eqref{eq:stab-MI+Ch} yields, for $t_n=nh\le T$,
	\begin{equation*}
		\|\varepsilon_n\|
		\le (1+Ch)^n\|\varepsilon_0\|
		+2C_{\rm loc}h^3\sum_{j=0}^{n-1}(1+Ch)^j.
	\end{equation*}
	Using $(1+Ch)^n\le e^{CT}$ for $nh\le T$ and
	\[
	\sum_{j=0}^{n-1}(1+Ch)^j=\frac{(1+Ch)^n-1}{Ch}\le \frac{(1+Ch)^n}{Ch}\le\frac{e^{CT}}{Ch},
	\]
	we obtain
	\[
	\|\varepsilon_n\|\le e^{CT}\|\varepsilon_0\|+\frac{2C_{\rm loc}}{C}e^{CT}h^2,
	\]
	which proves the stated estimate with $C_1=C$ and $C_2=2C_{\rm loc}/C$.
\end{proof}
Before presenting the second main theorem about the energy behaviour of EKahan method, we give the following lemma showing a property shared by a general exponential integrator based on the variation of constant formula applied to the semilinear system \eqref{semilinear HODE}. The proof of this lemma is inspired by Theorem 1 in \cite{Wuxinyuan2016}.
\begin{lemma}\label{lemma-energy}
	For  system \eqref{semilinear HODE}, any exponential integrator of the form 
	\begin{align*}
		x_{n+1}=e^{hA} x_{n} +h \phi(hA)Q\hat{\nabla}U(x_{n},x_{n+1}),
	\end{align*}
	satisfies the following  equation
	\begin{align*}
		\frac{1}{2}x_{n+1}^\mathsf{T} Mx_{n+1}-\frac{1}{2}x_{n}^\mathsf{T} Mx_{n}+(x_{n+1}-x_{n})^T\hat{\nabla}U(x_{n},x_{n+1})=0,
	\end{align*}
	where  $\hat{\nabla}U(x_{n},x_{n+1})$ is a discretization of the true gradient $\nabla U(x_{n})$.
\end{lemma}
{\color{black}Proofs of this lemma and the subsequent lemmas and corollaries are given in the Appendix. }

The discretized gradient $\hat{\nabla}U(x_{n},x_{n+1})$  include the discrete gradient obtained by AVF, the discretization given by Kahan's method as well as a more general discretization by Runge-Kutta method.  We denote the discretized gradient obtained via Kahan's method by $\nabla_{\mathrm{K}} U$, which can be written as 
\begin{align*}
	\nabla_{\mathrm{K}} U=-\frac{1}{2}\nabla U(x_n)+2\nabla U(\frac{x_n+x_{n+1}}{2})-\frac{1}{2}\nabla U(x_{n+1}).
\end{align*}
{\color{black}For the AVF choice, $\hat{\nabla}U$ is a discrete gradient and Lemma~\ref{lemma-energy} implies exact energy conservation.
	For the Kahan choice, Lemma~\ref{lemma-energy} will yield the following \emph{energy-increment} identity.}

\begin{theorem}\label{theorem-energy}
	For a Hamiltonian system \eqref{semilinear HODE} with a homogeneous cubic function $U(x)$, the EKahan scheme \eqref{EKahan scheme} satisfies the stepwise identity
	\begin{align}\label{step-wise-E}
		H_{n+1}-H_{n}=U(x_{n+1}-x_{n}),
	\end{align}
	where $H_{n}=\frac{1}{2}x_{n}^\mathsf{T} Mx_{n}+U(x_{n})$ is the discrete energy at $t=t_n$.
\end{theorem}
\begin{proof}
	Applying Lemma \ref{lemma-energy} with $\hat{\nabla}U=\nabla_{\mathrm K}U$ gives
	\begin{equation}\label{er-EH-Kahan}
		H_{n+1}-H_{n}=U(x_{n+1})-U(x_{n})-(x_{n+1}-x_{n})^{\mathsf{T}}\nabla_{\text{K}}U(x_{n},x_{n+1}).
	\end{equation}
	For a homogeneous cubic function $U(x)$, we have  $U(x)=\frac{1}{6}x^TU^{\prime\prime}(x)x$ and $\nabla_{\mathrm{K}} U=\frac{1}{2}U^{\prime\prime}(x_{n+1})x_{n}$. Then \eqref{er-EH-Kahan} can be reformulated as 
	\begin{align*}
		&H_{n+1}-H_{n}\\
		&=\frac{1}{6}{x_{n+1}}^\mathsf{T}U^{\prime\prime}(x_{n+1})x_{n+1}-\frac{1}{6}{x_{n}}^\mathsf{T}U^{\prime\prime}(x_{n})x_{n}-\frac{1}{2}{x_{n+1}}^\mathsf{T}U^{\prime\prime}(x_{n+1})x_{n}\\
		&+\frac{1}{2}{x_{n}}^\mathsf{T}U^{\prime\prime}(x_{n+1})x_{n}\\
		&=\frac{1}{6}{x_{n+1}}^\mathsf{T}U^{\prime\prime}(x_{n+1})x_{n+1}-\frac{1}{6}{x_{n+1}}^\mathsf{T}U^{\prime\prime}(x_{n+1})x_{n}+\frac{1}{6}{x_{n}}^\mathsf{T}U^{\prime\prime}(x_{n+1})x_{n}\\
		&-\frac{1}{6}{x_{n}}^\mathsf{T}U^{\prime\prime}(x_{n})x_{n}-\frac{1}{3}{x_{n+1}}^\mathsf{T}U^{\prime\prime}(x_{n+1})x_{n}+\frac{1}{3}{x_{n}}^\mathsf{T}U^{\prime\prime}(x_{n+1})x_{n}\\
		&=\frac{1}{6}(x_{n+1}-x_{n})^\mathsf{T}(U^{\prime\prime}(x_{n+1})x_{n+1}+U^{\prime\prime}(x_{n})x_{n}-2U^{\prime\prime}(x_{n+1})x_{n})\\
		&=\frac{1}{6}(x_{n+1}-x_{n})^\mathsf{T} U^{\prime\prime}(x_{n+1}-x_{n})(x_{n+1}-x_{n})\\
		&=U(x_{n+1}-x_{n}).
	\end{align*} 
\end{proof}
\begin{corollary}\label{cor:nonhom-cubic}
	{\color{black}Let the potential $U$ be a general cubic polynomial with homogeneous decomposition as
		\[
		U(x)=U_3(x)+U_2(x)+U_1(x)+U_0,
		\]
		where $U_j$ is homogeneous of degree $j$.
		Then the EKahan scheme \eqref{EKahan scheme} satisfies
		\[
		H_{n+1}-H_n = U_3(x_{n+1}-x_n),
		\]
		showing that the stepwise energy increment depends only on the cubic component.}
\end{corollary}
\begin{remark}\label{rem:energy-positioning}
	{\color{black}Theorem~\ref{theorem-energy} and Corollary~\ref{cor:nonhom-cubic} provide a computable description of the one-step energy defect in terms of the numerical increment. This complements standard finite-time estimates obtained indirectly from second-order convergence, and explains the favourable long-time behaviour of EKahan despite not implying exact energy conservation or symplecticity for general nonlinear Hamiltonian systems.
	}
\end{remark}
\subsection{\textbf{The multistep EKahan scheme for  higher-degree polynomial Hamiltonians}}\label{Method-continue}

{\color{black}We now extend the one-step EKahan scheme \eqref{EKahan scheme} from cubic potentials to polynomial potentials of degree $k+2$ with $k\in \mathbb{Z}^+$. 
	The key point is to retain the two features of EKahan: (i) exponential propagation of the linear flow, and (ii) a Kahan-type reversible discretization of the polynomial Hamiltonian part that requires solving only one linear system per time step. To achieve (ii), we portion the nonlinearity  by employing the symmetric $(k+1)$-linear polarization introduced in \cite{MR3452144}, leading to a $k$-step method.}

We denote the symmetric $(k+1)$-linear polarization of $\nabla U$ by 
\begin{equation}\label{function-polarization}
	\nabla_{\mathrm{K}} U(x_n,x_{n+1},\cdots,x_{n+k})
	=
	\sum_{\substack{1 \le m \le k+1 \\ n \le i_1 < \cdots < i_m \le n+k}}
	\frac{(-1)^{k+1-m}}{(k+1)!}\,\nabla U(x_{i_1}+\cdots+x_{i_m}).
\end{equation}
Then, we propose the following $k$-step EKahan scheme 
\begin{equation}\label{EKahan scheme-high-order}
	x_{n+k}=e^{khA} x_{n} +kh \phi(khA)Q\nabla_{\mathrm{K}} U(x_{n},x_{n+1},\cdots,x_{n+k}).
\end{equation}
{\color{black}\begin{lemma}\label{lemma-highorder-linearimplicit}
		Let $U$ be a polynomial of degree $k+2$. Then for each $n$ there exist a matrix $B_n\in\mathbb{R}^{d\times d}$ and a vector $\ell_n\in\mathbb{R}^d$, depending only on $(x_n,\dots,x_{n+k-1})$, such that
		\begin{equation}\label{eq:affine-polar}
			\nabla_{\mathrm{K}}U(x_n,\dots,x_{n+k}) = B_n x_{n+k} + \ell_n.
		\end{equation}
		Consequently, \eqref{EKahan scheme-high-order} is equivalent to the linear system
		\begin{equation*}
			\Bigl(I-kh\,\phi(khQM)\,Q B_n\Bigr)x_{n+k}
			= e^{khQM}x_n + kh\,\phi(khQM)\,Q\ell_n .
		\end{equation*}
		In particular, if $U$ is homogeneous, then $\ell_n=0$.
\end{lemma}}

{\color{black}We list several physically relevant polynomial potentials and give the corresponding  polarized gradient
	$\nabla_{\mathrm K}U$. }
\begin{itemize}
	\item For the standard lattice Duffing oscillator with on-site potential
	\[
	U(q)=\frac{\beta}{4}\sum_{j=1}^d (q^j)^4,
	\qquad \beta>0,
	\]
	the polarization formula \eqref{function-polarization} (with $k=2$) yields
	\[
	\nabla_{\mathrm K}U(q_n,q_{n+1},q_{n+2})
	=\beta\,(q_n\odot q_{n+1}\odot q_{n+2}),
	\]
	where $\odot$ denotes componentwise product, and the corresponding coefficient matrix and vector in \eqref{eq:affine-polar} are
	$B_n=\beta\,\mathrm{diag}(q_n\odot q_{n+1})$ and $\ell_n=0$.
	
	\item For the Fermi-Pasta-Ulam with quartic  potential
	\[
	U(q)=\frac{\beta}{4}\sum_{i=1}^{N-1}\bigl(r^i(q)\bigr)^4,
	\qquad r^i(q)=\frac{q^{i+1}-q^{i}}{\Delta x},
	\]
	introduce the strain vector $r(q)=(r^1(q),\ldots,r^{N-1}(q))^{\mathsf T}$ and $r_n=r(q_n)$.
	Then polarization formula \eqref{function-polarization} (with $k=2$) gives 
	\[
	\nabla_{\mathrm K}U(q_n,q_{n+1},q_{n+2})
	= D^{\mathsf T}\bigl(\beta\,(r_n\odot r_{n+1}\odot r_{n+2})\bigr),
	\]
	where $D$ is the forward difference operator, and the corresponding coefficient matrix and vector in \eqref{eq:affine-polar} are
	$B_n = D^{\mathsf T}\,\mathrm{diag}(w_n)\,D$
	and $\ell_n=0$, with $w_n=\beta\,(r_n\odot r_{n+1})$.
	\item For the H\'enon-Heiles potential
	\begin{equation*}
		U(p,q)=p^2 q-\frac{1}{3}q^3,
	\end{equation*}
	the two-point Kahan discretized gradient (with $k=1$) is
	\begin{equation*}
		\nabla_{\mathrm K}U(r_n,r_{n+1})
		=
		\begin{pmatrix}
			p_nq_{n+1}+p_{n+1}q_n\\
			p_np_{n+1}-q_nq_{n+1}
		\end{pmatrix}\qquad r=(p,q)^\mathsf{T},
	\end{equation*}
	and the corresponding coefficient matrix and vector in \eqref{eq:affine-polar} are
	$B_n=
	\begin{pmatrix}
		q_n & p_n\\
		p_n & -q_n
	\end{pmatrix}$ and $\ell_n=0$. 
\end{itemize}

{\color{black}Before presenting the properties of scheme \eqref{EKahan scheme-high-order}, we give the following notations and assumptions:
	\paragraph{Stacked state and history norm}
	Define the stacked state 
	\[
	Y_n = (x_{n+k-1},x_{n+k-2},\dots,x_n)\in(\mathbb{R}^d)^k,
	\]
	and the norm on the $k$-step history
	\[
	\|Y_n\|_\infty=\max_{0\le j\le k-1}\|x_{n+j}\|.
	\]
	\paragraph{$k$-step assumptions}
	For the $k$-step scheme~\eqref{EKahan scheme-high-order}, we assume that along the numerical trajectory
	the matrices $B_n$ in Lemma~\ref{lemma-highorder-linearimplicit} satisfy $\|B_n\|\le M_B$, and the uniform bound 
	$M_\phi=\max\limits_{0\le s\le kT}\|\phi(sA)\|$ holds for a fix final time $T>0$. 
	Finally, we assume the mappings $(x_n,\ldots,x_{n+k-1})\mapsto B_n$ and $(x_n,\ldots,x_{n+k-1})\mapsto \ell_n$
	are Lipschitz on the relevant bounded history set.}
{\color{black}\begin{theorem}\label{thm-highorder-wellposed}
		Assuming the $k$-step assumptions, then the $k$-step EKahan scheme \eqref{EKahan scheme-high-order} is well-defined for a sufficiently small $h$ and stable on $[0,T]$ with respect to starting perturbations.
\end{theorem}}
{\color{black}\begin{proof}
		By Lemma~\ref{lemma-highorder-linearimplicit}, \eqref{EKahan scheme-high-order} is equivalent to
		\begin{equation}\label{eq:linear-system-recall}
			\bigl(I-kh\,\phi(khQM)\,QB_n\bigr)x_{n+k}
			= e^{khQM}x_n + kh\,\phi(khQM)\,Q\ell_n .
		\end{equation}
		Under the $k$-step assumptions, there will exist a $h_0$ such that 
		\begin{equation}\label{eq:stepsize-sovl}
			kh\,M_\phi\,\|Q\|\,M_B \le \theta<1
			\quad\text{for all }0<h\le h_0.
		\end{equation}
		Then the coefficient matrix $I-kh\,\phi(khA)\,Q B_n$ is invertible and Neumann-series bound further gives 
		\begin{equation}\label{eq:inv-unif}
			\bigl\|\bigl(I-kh\,\phi(khA)\,Q B_n\bigr)^{-1}\bigr\|
			\le \frac{1}{1-\theta}
			\quad\text{for all }n\text{ and }0<h\le h_0 .
		\end{equation}
		Therefore $x_{n+k}$ is uniquely determined.
		
		Let $\{x_n\}$ and $\{\tilde x_n\}$ be two numerical solutions produced by \eqref{EKahan scheme-high-order} with the same step size. Set the error of the current state by $\delta_n=x_n-\tilde x_n$ and the error of the stacked state by $\Delta_n=Y_n-\tilde Y_n$ .
		Define
		\[
		G_n=kh\,\phi(khQM)\,QB_n,\qquad
		b_n=e^{khQM}x_n+kh\,\phi(khQM)\,Q\ell_n,
		\]
		and $\tilde G_n,\tilde b_n$ analogously. Then from \eqref{eq:linear-system-recall} we have 
		$(I-G_n)x_{n+k}=b_n$ and $(I-\tilde G_n)\tilde x_{n+k}=\tilde b_n$.
		Subtracting and using the resolvent identity yields
		\begin{equation}\label{eq:delta-resolvent}
			\delta_{n+k}
			= (I-G_n)^{-1}(b_n-\tilde b_n)
			+ (I-G_n)^{-1}(G_n-\tilde G_n)(I-\tilde G_n)^{-1}\tilde b_n .
		\end{equation}
		Taking $h_0$ such that $\theta=\frac{1}{2}$ in \eqref{eq:stepsize-sovl}-\eqref{eq:inv-unif}, we have $ kh\,M_\phi\,\|Q\|\,M_B \le \tfrac12$
		and 
		\begin{equation}\label{eq:inv-1+Ch}
			\|(I-G_n)^{-1}\|\le \frac{1}{1-\|G_n\|}\le 1+2\|G_n\|\le 1+C_B h,
			\quad C_B:=2kM_\phi\|Q\|M_B,
		\end{equation}
		and the same bound holds for $\|(I-\tilde G_n)^{-1}\|$.
		
		Moreover, recalling from the proof of Theorem \ref{thm:global-error}, we have 
		$\|e^{khQM}\|\le 1+C_E h$ with $C_E=k\|QM\|M_E$.
		Taking norms in \eqref{eq:delta-resolvent} and using \eqref{eq:inv-1+Ch} yield
		\begin{equation}\label{eq:delta-bound1}
			\|\delta_{n+k}\|
			\le (1+C_B h)\,\|b_n-\tilde b_n\|
			+ (1+C_B h)^2\,\|G_n-\tilde G_n\|\,\|\tilde b_n\|.
		\end{equation}
		Furthermore, we have
		\begin{equation*}
			\begin{split}
				\|b_n-\tilde b_n\|
				&\le \|e^{khQM}\|\,\|\delta_n\| + kh\,\|\phi(khQM)\|\,\|Q\|\,\|\ell_n-\tilde\ell_n\|\\
				&\le (1+C_E h)\|\delta_n\| + kh\,M_\phi\|Q\|\,\|\ell_n-\tilde\ell_n\|,
			\end{split}
		\end{equation*}
		and
		\[
		\|G_n-\tilde G_n\|
		\le kh\,\|\phi(khQM)\|\,\|Q\|\,\|B_n-\tilde B_n\|
		\le kh\,M_\phi\|Q\|\,\|B_n-\tilde B_n\|.
		\]
		By the Lipschitz assumption on the mappings $(x_n,\ldots,x_{n+k-1})\mapsto B_n$ and
		$(x_n,\ldots,x_{n+k-1})\mapsto \ell_n$, there exist constants $L_B,L_\ell>0$ such that
		\[
		\|B_n-\tilde B_n\|\le L_B\|\Delta_n\|_\infty,\qquad
		\|\ell_n-\tilde\ell_n\|\le L_\ell\|\Delta_n\|_\infty.
		\]
		In addition, boundedness of the trajectories implies $\|\tilde b_n\|\le M_b$ for some constant $M_b>0$. Substituting these estimates into \eqref{eq:delta-bound1} and using $\|\delta_n\|\le\|\Delta_n\|_\infty$ gives
		\[
		\|\delta_{n+k}\|
		\le \Bigl((1+C_B h)(1+C_E h+khM_\phi\|Q\|L_\ell)
		+ kh(1+C_B h)^2M_\phi\|Q\|L_B\,M_b\Bigr)\|\Delta_n\|_\infty.
		\]
		Therefore, for $0<h\le h_0$ there exists $\widetilde C>0$ independent of $h$ such that
		\[
		\|\delta_{n+k}\|\le (1+\widetilde C h)\|\Delta_n\|_\infty.
		\]
		Consequently,
		\begin{equation}\label{high-esti-one-vec}
			\|\Delta_{n+1}\|_\infty
			=\max\bigl\{\|\delta_{n+k}\|,\ \|\delta_{n+k-1}\|,\ldots,\|\delta_{n+1}\|\bigr\}
			\le (1+\widetilde C h)\|\Delta_n\|_\infty.
		\end{equation}
		Iterating and using $nh\le T$ yields
		\[
		\|\Delta_n\|_\infty \le (1+\widetilde C h)^n\|\Delta_0\|_\infty \le e^{\widetilde C T}\|\Delta_0\|_\infty,
		\qquad 0\le nh\le T,
		\]
		which completes the proof.
\end{proof}}
\medskip
Moreover, the properties of EKahan scheme established for systems with a cubic Hamiltonian can be extended to higher-degree polynomial Hamiltonians, as summarized next.

\begin{theorem}\label{corollary-high-order-symmetric-order}
	For semilinear system \eqref{semilinear HODE} with polynomial potential $U(x)$ of degree $k+2$, the EKahan scheme \eqref{EKahan scheme-high-order} is symmetric.
\end{theorem}
\begin{proof}
	From equation \eqref{EKahan scheme-high-order} we get
	\begin{equation}\label{eq:highorder-map}
		x_{n+k}-e^{khA}x_n
		=kh\,\phi(khA)\,Q\,\nabla_{\mathrm K}U(x_n,\ldots,x_{n+k}),
		\qquad A=QM .
	\end{equation}
	Using $\phi(-Z)=e^{-Z}\phi(Z)$ and $e^{-khA}=(e^{khA})^{-1}$, we can rewrite \eqref{eq:highorder-map} as
	\begin{equation*}
		x_n-e^{-khA}x_{n+k}
		=-kh\,\phi(-khA)\,Q\,\nabla_{\mathrm K}U(x_n,\ldots,x_{n+k}).
	\end{equation*}
	Since $\nabla_{\mathrm K}U$ is symmetric in its arguments, $\nabla_{\mathrm K}U$ will not change if we reverse the sequence
	$(x_n,\ldots,x_{n+k})$.
	Therefore, the method is symmetric.
\end{proof}
{\color{black}
	\begin{lemma}\label{lem:exp-mid-quad}
		Let $A\in\mathbb R^{d\times d}$ and $g\in C^2([0,kh];\mathbb R^d)$. For all $t\in[0,kh]$, we claim
		\[
		\int_0^{kh} e^{(kh-s)A}g(s)\mathrm{d}s
		=kh\,\phi(khA)\,g(kh/2)+R,
		\qquad \|R\|\le C h^3,
		\]
		where $C$ depends on $k$, $M_E$, $\|A\|$, and $\max\limits_{0\le s\le kh}\bigl(\|g(s)\|+\|g'(s)\|+\|g''(s)\|\bigr)$.
	\end{lemma}
	\begin{lemma}\label{lem:polar-mid}
		Assume $U$ is homogeneous and let $x(t)$ be a smooth solution on $[t_n,t_n+kh]$.
		Set $\tilde x_{n+i}:=x(t_n+ih)$ and $t_{n+\frac{k}{2}}:=t_n+\frac{k}{2}h$.
		Then there exists $C>0$, independent of $h$, such that
		\[
		\bigl\|\nabla_{\mathrm K}U(\tilde x_n,\ldots,\tilde x_{n+k})
		-\nabla U\bigl(x(t_{n+\frac{k}{2}})\bigr)\bigr\|
		\le C h^2 .
		\]
	\end{lemma}
	\begin{lemma}\label{lem:local-defect-high}
		Let $x(t)$ solve \eqref{semilinear HODE} with $U$ homogeneous of degree $k+2$ and set
		$\tilde x_{n+i}=x(t_n+ih)$. Then the local defect of \eqref{EKahan scheme-high-order} satisfies
		\[
		\|\tau_{n+k}\|\le C_{\mathrm{loc}}\,h^3,
		\qquad t_n\le T,
		\]
		for a constant $C_{\mathrm{loc}}$ independent of $h$.
	\end{lemma}
	\begin{proof}
		By the variation of constants formula for $\dot x=Ax+Q\nabla U(x)$, we have 
		\[
		\tilde x_{n+k}=e^{khA}\tilde x_n+\int_{0}^{kh}e^{(kh-s)A}Q\nabla U(x(t_n+s))\mathrm{d}s.
		\]
		Hence
		\[
		\tau_{n+k}
		=\int_{0}^{kh}e^{(kh-s)A}Q\nabla U(x(t_n+s))\mathrm{d} s
		-kh\,\phi(khA)\,Q\,\nabla_{\mathrm K}U(\tilde x_n,\ldots,\tilde x_{n+k}).
		\]
		Add and subtract the exponential midpoint quantity $kh\,\phi(khA)\,Q\,\nabla U(x(t_{n+\frac{k}{2}}))$:
		\begin{equation*}
			\begin{split}
				\tau_{n+k}&= \underbrace{\Bigl[\int_{0}^{kh}e^{(kh-s)A}Q\nabla U(x(t_n+s))\mathrm{d} s
					-kh\,\phi(khA)\,Q\,\nabla U(x(t_{n+\frac{k}{2}}))\Bigr]}_{=:I_1}\\
				&+\underbrace{kh\,\phi(khA)\,Q\Bigl[\nabla U(x(t_{n+\frac{k}{2}}))
					-\nabla_{\mathrm K}U(\tilde x_n,\ldots,\tilde x_{n+k})\Bigr]}_{=:I_2}.
			\end{split}
		\end{equation*}
		For $I_1$, apply Lemma~\ref{lem:exp-mid-quad} with $g(s)=Q\nabla U(x(t_n+s))$ to get $\|I_1\|\le C h^3$.
		For $I_2$, use the boundedness of $\phi(khA)$ on $[0,T]$ and Lemma~\ref{lem:polar-mid} to obtain
		\[
		\|I_2\|\le kh\,\|\phi(khA)\|\,\|Q\|\,
		\bigl\|\nabla_{\mathrm K}U(\tilde x_n,\ldots,\tilde x_{n+k})
		-\nabla U(x(t_{n+\frac{k}{2}}))\bigr\|
		\le C\,h\cdot h^2 = C h^3.
		\]
		Combining the two bounds yields $\|\tau_{n+k}\|\le C_{\mathrm{loc}} h^3$.
	\end{proof}
	\begin{lemma}\label{thm:zero-stab-high}
		Assume the exact solution stays in a bounded set $\mathcal D$ on $[0,T]$ and $U$ is polynomial.
		Then there exist $h_0>0$ and $C>0$ such that for all $0<h\le h_0$ the $k$-step EKahan update
		defines a unique one-step map $Y_{n+1}=\Phi_h(Y_n)$ on $\mathcal D^k$, and
		\[
		\|\Phi_h(Y)-\Phi_h(Z)\|_\infty\le (1+Ch)\,\|Y-Z\|_\infty,
		\qquad \forall\,Y,Z\in\mathcal D^k,
		\]
	\end{lemma}
	\begin{theorem}\label{thm:global-high}
		Under the assumptions of Lemma \ref{thm:zero-stab-high}, suppose the starting values
		$x_0,\ldots,x_{k-1}$ satisfy $\|x_j-\tilde x_j\|\le C h^2$ for $j=0,\ldots,k-1$.
		Then for $t_n=nh\le T$ and $h\le h_0$,
		\[
		\max_{0\le n\le N}\|x_n-\tilde x_n\|\le C_T\,h^2,
		\]
		i.e., the $k$-step EKahan scheme is globally second order on finite time intervals.
	\end{theorem}
	\begin{proof}
		Inserting the exact solution values $\tilde Y_n$ into the scheme gives
		\begin{equation*}
			\tilde Y_{n+1}=\Phi_h(\tilde Y_n)+E_{n+1},
			\qquad \|E_{n+1}\|\le C_{\mathrm{loc}}h^3.
		\end{equation*}
		where $E_n$ contains only the last component defect $\tau_{n+k}$ from Lemma~\ref{lem:local-defect-high}. 
		Subtracting from $Y_{n+1}=\Phi_h(Y_n)$ and using Lemma~\ref{thm:zero-stab-high} yields the following  recursion  for the  extended error $\Delta_n=Y_n-\tilde Y_n$
			\begin{equation*}
				\|\Delta_{n+1}\|_\infty \le (1+Ch)\|\Delta_n\|_\infty+C_{\mathrm{loc}}h^3.
			\end{equation*}
		Applying the discrete Gronwall inequality gives, for $t_n\le T$,
		\[
		\|\Delta_n\|_\infty \le e^{CT}\|\Delta_0\|_\infty + \frac{C_{\mathrm{loc}}}{C}\,(e^{CT}-1)\,h^2.
		\]
		Since $\|\Delta_0\|_\infty=\max\limits_{0\le j\le k-1}\|x_j-\tilde x_j\|=O(h^2)$ by assumption, the claim follows.
\end{proof}}
\begin{corollary}\label{corollary-high-order-energy}
	For the semilinear system \eqref{semilinear HODE}, any exponential integrator of the form
	\[
	x_{n+k}=e^{khA} x_{n} +kh \phi(khA)\,Q\,\hat{\nabla}U(x_{n},x_{n+1},\cdots,x_{n+k}),
	\]
	where  $\hat{\nabla}U(x_{n},x_{n+1},\cdots,x_{n+k})$ is a discretization of the true gradient $\nabla U(x_n)$, satisfies
	\[
	\frac{1}{2k}x_{n+k}^\mathsf{T} Mx_{n+k}-\frac{1}{2k}x_{n}^\mathsf{T} Mx_{n}
	+\frac{1}{k}(x_{n+k}-x_{n})^\mathsf{T} \hat{\nabla}U(x_{n},x_{n+1},\cdots,x_{n+k})=0.
	\]
\end{corollary}
\begin{corollary}\label{corollary-EKahan-energy-high-order}
	For a homogeneous function $U(x)$ of degree $k+2$, EKahan scheme \eqref{EKahan scheme-high-order} satisfies the stepwise identity
	\begin{equation}\label{higher-degree-EI}
		\begin{split}
			H_{n+1}-H_{n}=&\bar{U}(x_{n+1},x_{n+2},\dots,x_{n+k-1},x_{n+k},x_{n+k},x_{n+1}-x_{n})\\
			&+\bar{U}(x_{n},x_{n+1},\dots,x_{n+k-2},x_{n+k-1},x_{n},x_{n+k}-x_{n+k-1})\\
			&-2\bar{U}\Bigl(x_{n},x_{n+1},\dots,x_{n+k-1},x_{n+k},\frac{x_{n+k}-x_{n}}{k}\Bigr),
		\end{split}
	\end{equation}
	where
	\[
	H_{n}=\frac{1}{2k}\sum_{i=0}^{k-1} x_{n+i}^\mathsf{T} Mx_{n+i}
	+\bar{U}(x_{n},x_{n+1},\dots,x_{n+k-2},x_{n+k-1},x_{n},x_{n+k-1})
	\]
	is the discrete energy at $t=t_n$, and $\bar{U}(\cdot,\dots,\cdot)$ is a symmetric $(k+2)$-linear form associated with $U$ satisfying $\bar{U}(x,\dots,x)=U(x)$.
\end{corollary}
\begin{remark}
	For a nonhomogeneous polynomial $U(x)$, we homogenize it by introducing an auxiliary variable $x^0$ and an extended state $\bar x=(x^0,x^1,\ldots,x^m)^\mathsf{T}$ so that $\bar U(1,x)=U(x)$. To preserve the structural properties of the original system, we also augment  $Q$ to $\bar{Q}$ and $M$ to $\bar{M}$ by adding a zero row and column. Notably, $\bar{Q}$ remains skew-symmetric and $\bar{M}$ remains symmetric. Furthermore, applying the EKahan method to 
	\[
	\dot{\bar{x}}(t) = \bar{Q} \bar{M} \bar{x}(t) + \bar{Q} \nabla \bar{U}(\bar{x}(t)), \quad \bar{x}(t_0) = [1, x_0^\mathsf{T}]^\mathsf{T},
	\]
	produces the same update for $x$ as applying the EKahan method directly to the original system \eqref{semilinear HODE}. Therefore, all results previously established for systems with homogeneous polynomial Hamiltonians  extend naturally to those involving nonhomogeneous polynomials.
\end{remark}

\section{Numerical test}\label{Experiment}

In this section, we apply the proposed EKahan scheme to  three systems including an ODE and two PDEs with high dimension and high-order nonlinear functions. We compare EKahan scheme  with the state-of-the-art symmetric structure-preserving methods including the following three schemes:
\begin{itemize}
	\item Kahan: the Kahan method  presented in \cite{celledoni2012geometric};
	\item EAVF: the exponential Averaged Vector Field method presented in \cite{Wuxinyuan2016};
	\item LIEEP: the linearly implicit energy-preserving exponential scheme illustrated in \cite{li2022new}.
\end{itemize}
We define the global error at time $ t_n = t_0 + nh $ as
\begin{equation*}
	\mathcal{E}_G = \max_{n \geq 0} \Vert y_{n} - y(t_{n}) \Vert,
\end{equation*}
where \( y(t_{n}) \) is the reference exact solution. For EKahan, we measure the residual of the stepwise energy-increment identity by
\begin{equation*}\label{Residualerror-plot}
	e_n=H_{n+1}-H_n-G_n,
\end{equation*}
where $G_n=U(x_{n+1}-x_n)$ for cubic Hamiltonians and  $ H_n $ is the discrete energy defined in Theorem \ref{theorem-energy}. For higher-degree polynomial Hamiltonians, $G_n$ is given by the right-hand side of equation~\eqref{higher-degree-EI} and $ H_n $ is defined as in Corollary~\ref{corollary-EKahan-energy-high-order}.
The energy error  is defined as
\begin{equation}\label{Eerror-plot}
	\mathcal{E}_H(n) = \lvert H_{n} - H(y(t_n)) \rvert.
\end{equation}
The reference exact solution is obtained using the 6-order continuous Runge-Kutta (CRK) method \cite{hairer2010energy} with the form
\begin{equation*}
	\left\{\begin{array}{l}
		y_{n+1 / 3}=y_n+h Q \int_0^1\left(\frac{37}{27}-\frac{32}{9} \sigma+\frac{20}{9} \sigma^2\right) \nabla H\left(Y_\sigma\right) \mathrm{d} \sigma \\
		y_{n+2 / 3}=y_n+h Q \int_0^1\left(\frac{26}{27}+\frac{8}{9} \sigma-\frac{20}{9} \sigma^2\right) \nabla H\left(Y_\sigma\right) \mathrm{d} \sigma \\
		y_{n+1}=y_n+h Q \int_0^1 \nabla H\left(Y_\sigma\right) \mathrm{d} \sigma
	\end{array}\right.,
\end{equation*}
where
\begin{align*}
	Y_\sigma= & -\frac{(3 \sigma - 1)(3 \sigma-2)(\sigma-1)}{2} y_n+\frac{3 \sigma(3 \sigma-2)(3 \sigma-3)}{2} y_{n+1 / 3} \\
	& -\frac{3 \sigma(3 \sigma-1)(3 \sigma-3)}{2} y_{n+2 / 3}+\frac{\sigma(3 \sigma-1)(3 \sigma-2)}{2} y_{n+1},
\end{align*}
and the integrals are evaluated exactly by the 5-point GL quadrature.  For both LIEEP method and EKahan scheme \eqref{EKahan scheme-high-order}, multiple initial starting points are required. These starting values are obtained using the aforementioned reference exact solution. The matrix exponential functions $e^{khQM}$ and $\phi(khQM)$ are computed using the MATLAB package presented in \cite{1206044} when no explicit form is given, where Pad\'e approximations are employed. {\color{black}We solve the linear systems arising from each step of the linearly implicit schemes (EKahan,  Kahan and LIEEP) by MATLAB's backslash operator. The nonlinear updates required for the reference solutions and the EAVF numerical solutions are computed via a fixed-point iteration, terminated when the 2-norm of the difference between successive iterates falls below a tolerance of $10^{-14}$. We also tested tolerances from $10^{-14}$ to $10^{-6}$ for EAVF and observed that EKahan is consistently the most efficient method across our test problems, except for the FPU system with $p=2$, where EKahan is not as efficient as EAVF when a bigger tolerance is considered.}  All numerical experiments were conducted using MATLAB (R2022b) on a 2020 MacBook Air equipped with an Apple M1 processor and 16 GB of unified memory.
\paragraph{\textbf{H\'enon-Heiles system}}
The H\'enon-Heiles model \cite{MR1128990} has been used to describe stellar motion in a long time course of gravitational potential energy tracking of galaxies with columnar symmetry. It has the following form
\begin{align*}\label{HHH}
	\left\{\begin{array}{ll}
		{\ddot q_1 =-q_1 -2 Dq_1 q_2},  \\ 
		{\ddot q_{2}=- q_2 -D q_1^{2}+C q_2^{2}},  \\ 
	\end{array}\right.
\end{align*}
where the independent variable $t \in [0,T]$. We consider $C=1$ and $D=1$. By introducing the momentum variables $\dot{q}_1 =p_1 ,\; \dot{q}_2 =p_2$, the system can be rewritten in the Hamiltonian form 
\begin{equation*}
	\dot{y}=Q \nabla H, \quad
	H(y) = \frac{1}{2} y^{\mathsf{T}} M y + U(y),
\end{equation*}
where \( y = (q_1, q_2, p_1, p_2)^{\mathsf{T}} \), \( Q \) is the canonical skew-symmetric matrix, \( M = I_4 \) denotes the \( 4 \times 4 \) identity matrix, and the potential function is defined as \( U(y) = q_1^2 q_2 - \frac{1}{3} q_2^3 \).  The exponential function can be computed explicitly by 
\begin{equation*}
	e^{hQM}= \begin{pmatrix}
		\cos(h)I_{2} & \sin(h)I_{2}\\
		-\sin(h)I_{2} & \cos(h)I_{2} 
	\end{pmatrix}.
\end{equation*}



\begin{figure}[H]
	\centering
	\begin{subfigure}[b]{0.4\textwidth}
		\includegraphics[width=\textwidth]{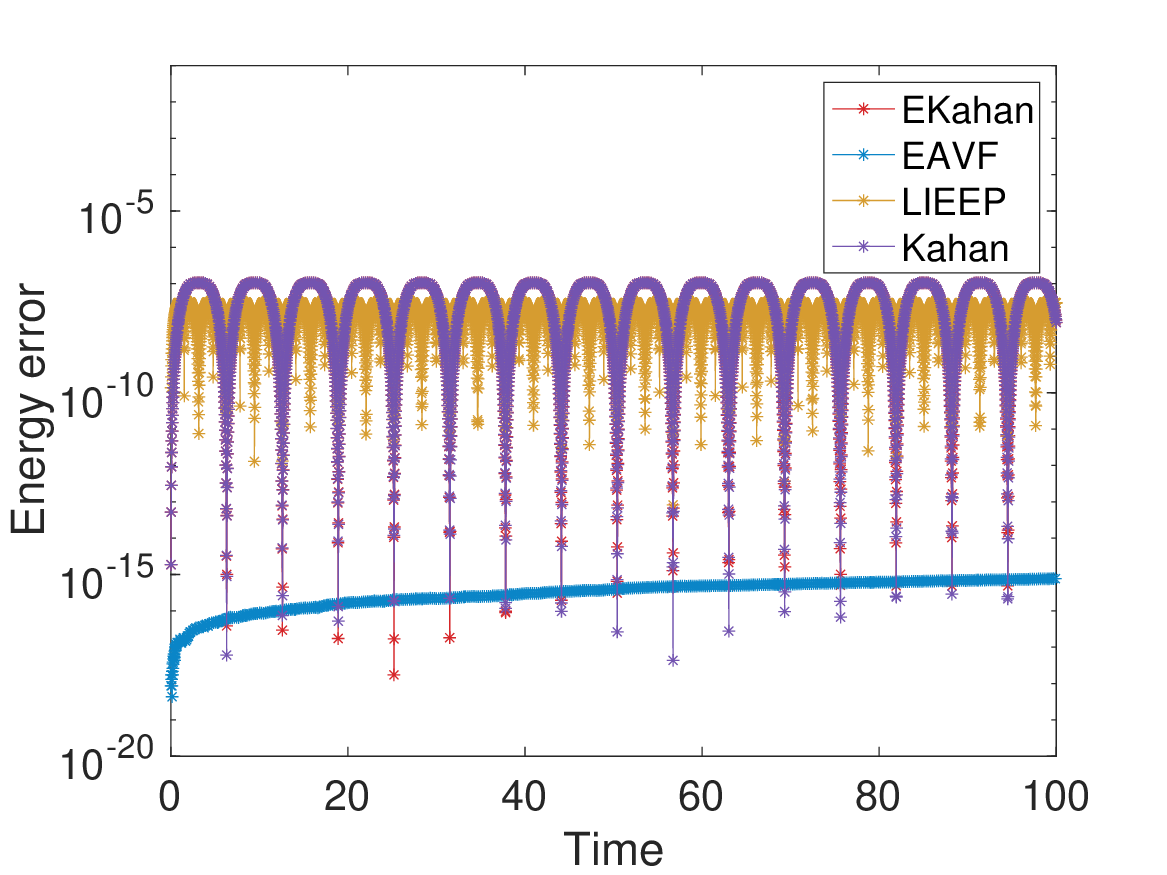}
		\caption[a]{}
		\label{fig:HH energy error}
	\end{subfigure}
	\begin{subfigure}[b]{0.4\textwidth}
		\includegraphics[width=\textwidth]{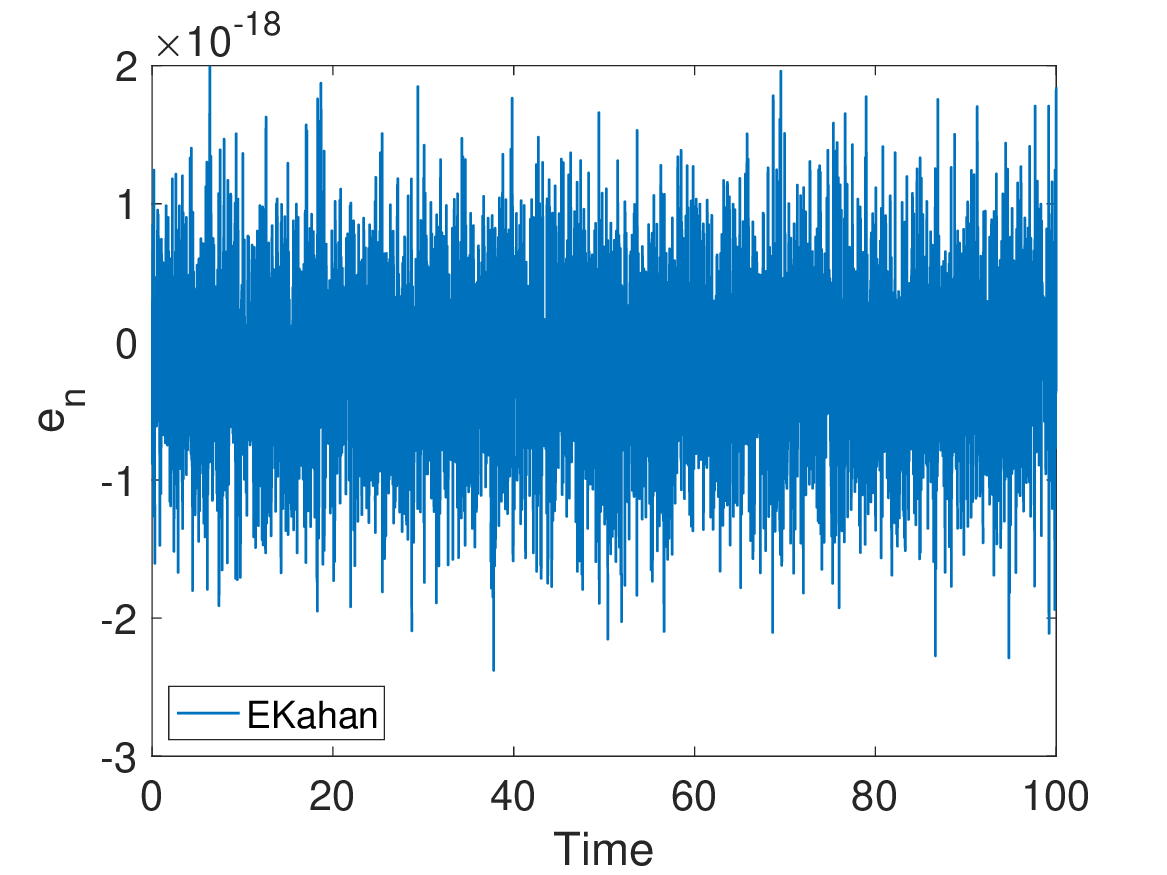}
		\caption[a]{}
		\label{fig:HH energy est}
	\end{subfigure}
	\caption{H\'enon-Heiles equation with $T=100$ and $h_{0}=0.02$. Left: energy error defined in \eqref{Eerror-plot}; right:  the residual of equation \eqref{step-wise-E}  in Theorem \ref{theorem-energy}.}
	\label{fig:HH solution energy}
\end{figure}

Consider the initial data $y_{0}=(0,-0.082,0,0)^\mathsf{T}$, $T=100$ and time step size $h_{i}=0.02/2^{i},\;i=0,1,\dots,4$. Figure \ref{fig:HH energy error} demonstrates that the EAVF method preserves the discrete energy exactly and the energy errors of all linearly implicit methods remain bounded, varying from $10^{-18}$ to $10^{-7}$. Figure \ref{fig:HH energy est} verifies the result about the energy deviation at each time step by EKahan  shown in Theorem \ref{theorem-energy}. Furthermore, we verified numerically that the per-step energy increment of EKahan scales as $O(h^3)$, consistent with the stepwise identity in Theorem~\ref{theorem-energy}, and we omitted the corresponding figure due to space limitations.

\begin{figure}[H]
	\centering
	\begin{subfigure}[b]{0.4\textwidth}
		\includegraphics[width=\textwidth]{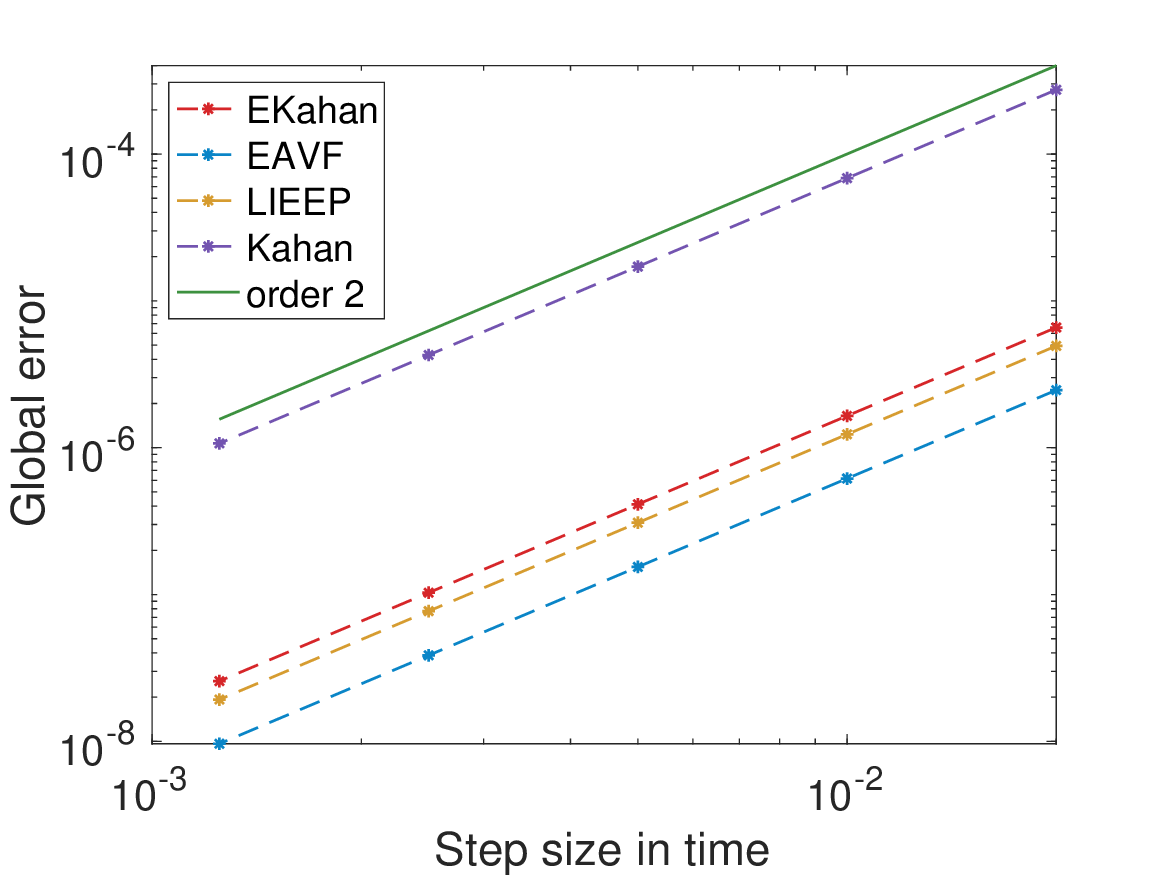}
		\caption[a]{order plot}
		\label{fig:HH order}
	\end{subfigure}
	\begin{subfigure}[b]{0.4\textwidth}
		\includegraphics[width=\textwidth]{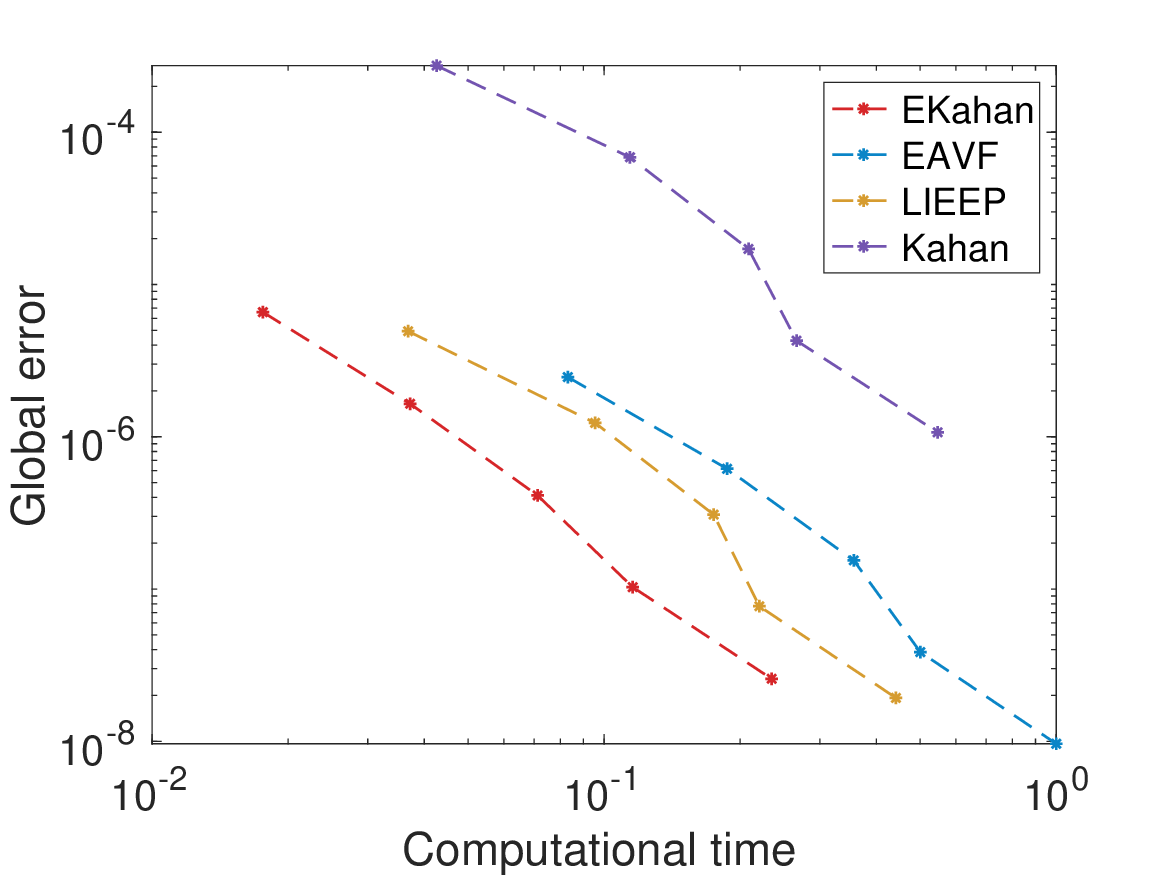}
		\caption[a]{efficiency plot}
		\label{fig:HH efficiency}
	\end{subfigure}   
	\caption{Global error and computational time of different methods for the H\'enon-Heiles equation with $T=100$ and different time step sizes. }
	\label{fig:HH order efficiency}
\end{figure}

In Figure \ref{fig:HH order efficiency}, we consider the global errors at $T=100$ and the computational cost using different step sizes. Figure \ref{fig:HH order} confirms that the EKahan method is second-order, and Figure \ref{fig:HH efficiency} shows that EKahan method is more efficient compared to the fully implicit EAVF method and the other two types of linearly implicit methods. In addition, we observe from Figure \ref{fig:HH order} that the exponential integrators are more accurate than Kahan's method since the linear part are exactly integrated using exponential integrators.

\paragraph{\textbf{FPU system with a cubic Hamiltonian}}
The $\alpha$-Fermi-Pasta-Ulam (FPU) system, originating from the study by Fermi, Pasta, Ulam, and Tsingou, has evolved into a paradigmatic model in nonlinear wave theory \cite{4376203}. We study a continuum generalization of the FPU system
\begin{equation*}\label{fpu1}
	\frac{\partial^2 u}{\partial t^2}=\beta \frac{\partial^3 u}{\partial t \partial^2 x} + \frac{\partial^2 u}{\partial x^2}[1+\epsilon (\frac{\partial u}{\partial x})^p] -\gamma \frac{\partial u}{\partial t} -m^2 u,
\end{equation*}
where the independent variables $(x,t) \in [0,L] \times [0,T]$, the parameter $\epsilon > 0$, the internal damping coefficient $\beta \geq0$, and the external damping coefficient $\gamma \geq0$. By introducing $\frac{\partial u}{\partial t}=v$ and  
setting $y=[u^\mathsf{T},v^\mathsf{T}]^{\mathsf{T}}$, we can rewrite the above equation  into the following Hamiltonian form 
\begin{equation*}
	\frac{\partial y}{\partial t} = \mathcal{Q}\frac{\delta \mathcal{H}}{\delta y},
\end{equation*}
where
\begin{equation*}
	\mathcal{Q}=
	\left[
	\begin{array}{cc}
		0 & 1 \\
		-1 & \beta\partial^2_x -\gamma \\
	\end{array}
	\right],
	\quad
	\mathcal{H}=\int_{0}^{L} H(t,u,v,u_x) \mathrm{d}x,
\end{equation*}
and the local energy density $H(t,u,v,u_x)$  is given by
\begin{equation*}
	H(t,u,v,u_x)=\frac{1}{2}u_x^{2}+\frac{m^2}{2}u^2+\frac{1}{2}v^2+\epsilon\frac{u_x^{p+2}}{(p+1)(p+2)}.
\end{equation*}
We consider $p=1$ and  homogeneous Dirichlet boundary conditions $u(0,t) = u(L,t) = 0$. Employing the second-order central difference to discretize the second derivative $\partial_x^2$ and the forward difference to discretize $\partial_x$, we derive the semi-discrete system 
\begin{equation*}
	\dot{y}=Q(My+\nabla U(y)),
\end{equation*}
where
\begin{equation*}
	Q=\left[
	\begin{array}{cc}
		0 & I \\
		-I & \beta D-\gamma I \\
	\end{array}
	\right],
	M=\left[
	\begin{array}{cc}
		m^2I-D & 0 \\
		0 & I  \\ 
	\end{array}
	\right],
	U(y)=\sum_{j=0}^{N-1} \frac{\epsilon}{6}(\frac{u^{j+1}-u^j}{\Delta x})^3.
\end{equation*}

We set  $m=0$, $\epsilon=\frac{3}{4}$, $L=128$, $T=100$, the spatial discretization to be $\Delta x=1$ and the time step size to be $h_{i}=1/2^{i}$, $i=1,2,\ldots,4$.  
The initial conditions are given by
\begin{equation*}
	u^j(0) = q^j(0), \quad v^j(0) = \dot{q}^j(0),
\end{equation*}
where the function $q^j(t)$ takes the form
\begin{equation*}
	q^j(t) = 5\sum_{k\in\{32,96\}} \ln \frac{1+e^{2(\alpha(j-k)+t\sinh\alpha)}}{1+e^{2(\alpha(j-k-1)+t\sinh\alpha)}}
\end{equation*}
with $\alpha=0.1$.

\begin{figure}[H]
	\centering
	\begin{subfigure}[b]{0.4\textwidth}
		\includegraphics[width=\textwidth]{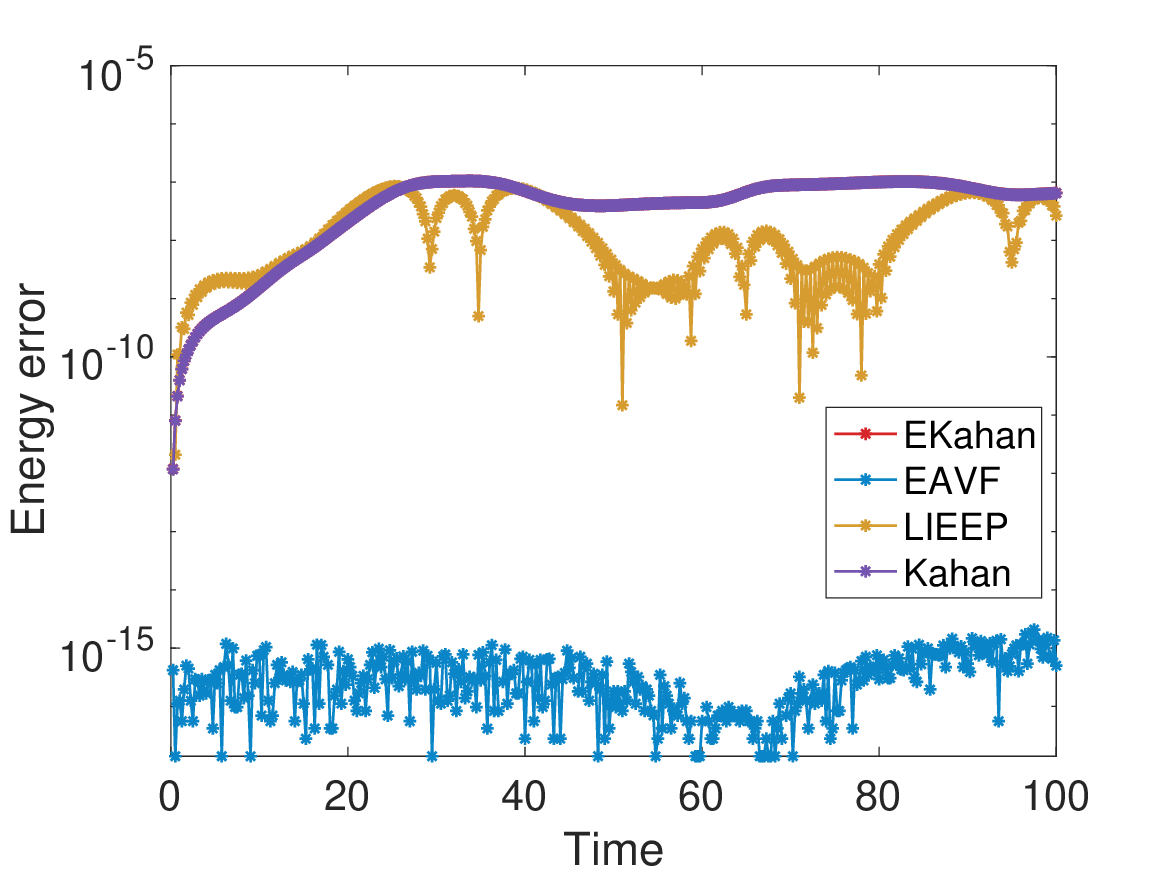}
		\caption[a]{$\gamma=0$, $\beta=0$}
		\label{fig:FPU energy_error_gamma=0_beta=0}
	\end{subfigure}
	\begin{subfigure}[b]{0.4\textwidth}
		\includegraphics[width=\linewidth]{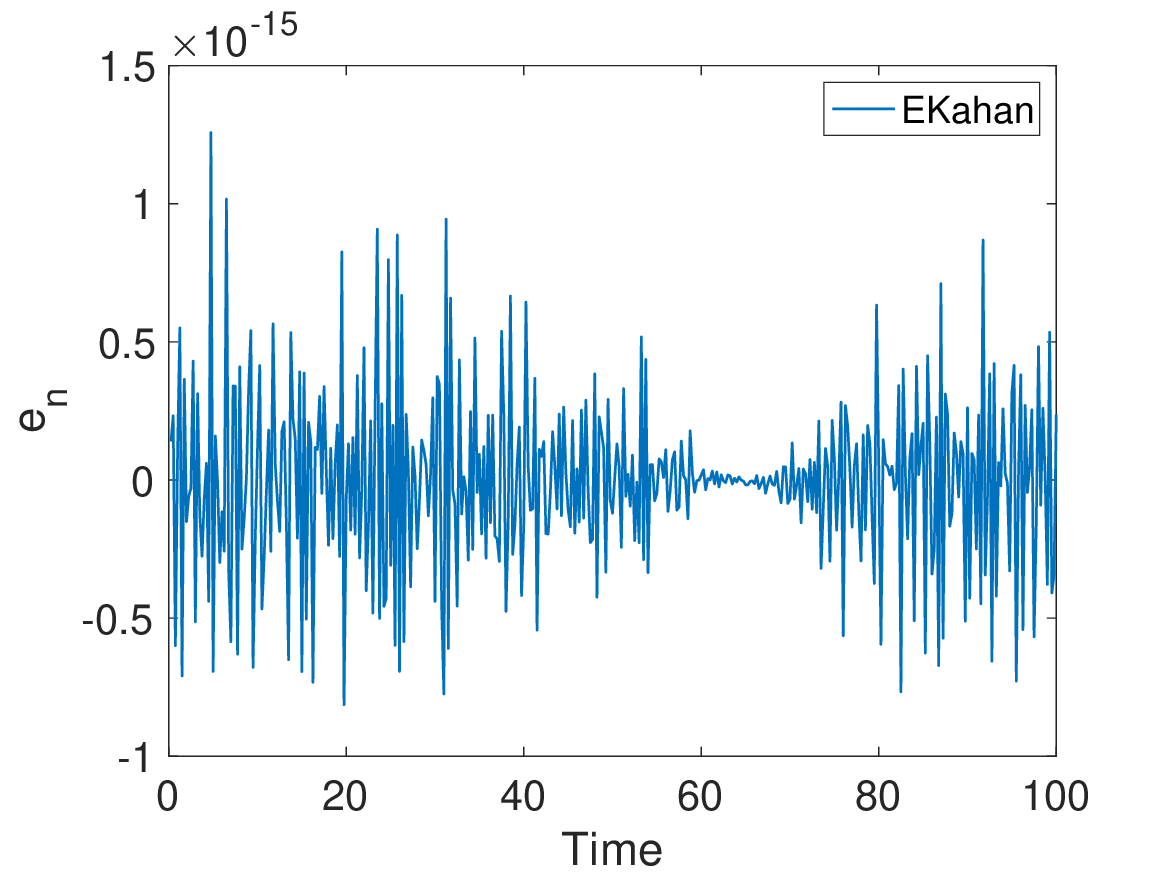}
		\caption[b]{$\gamma=0$, $\beta=0$}
		\label{fig:FPU est}
	\end{subfigure}
	\begin{subfigure}[b]{0.4\textwidth}
		\includegraphics[width=\textwidth]{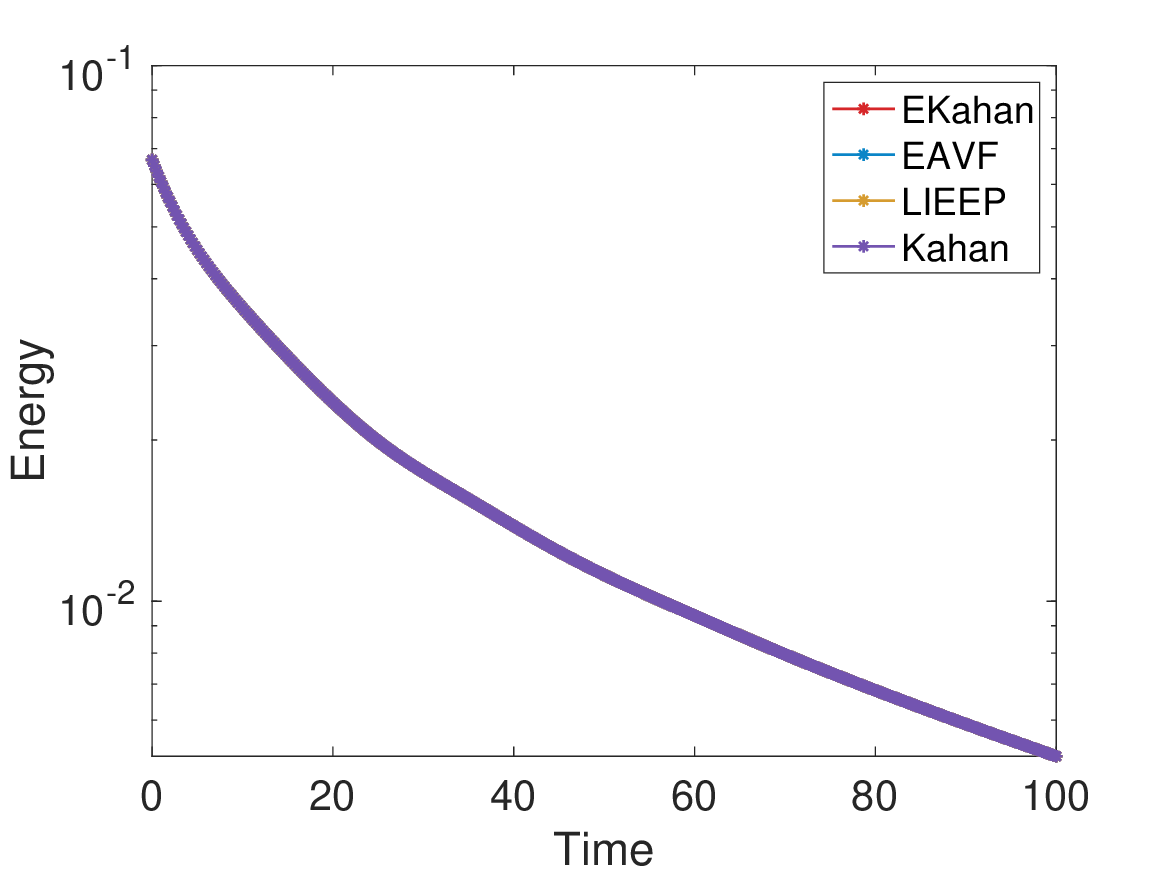}
		\caption[a]{$\gamma=0.1$, $\beta=0$}
		\label{fig:FPU energy_error_gamma=0.1_beta=0}
	\end{subfigure}
	\begin{subfigure}[b]{0.4\textwidth}
		\includegraphics[width=\textwidth]{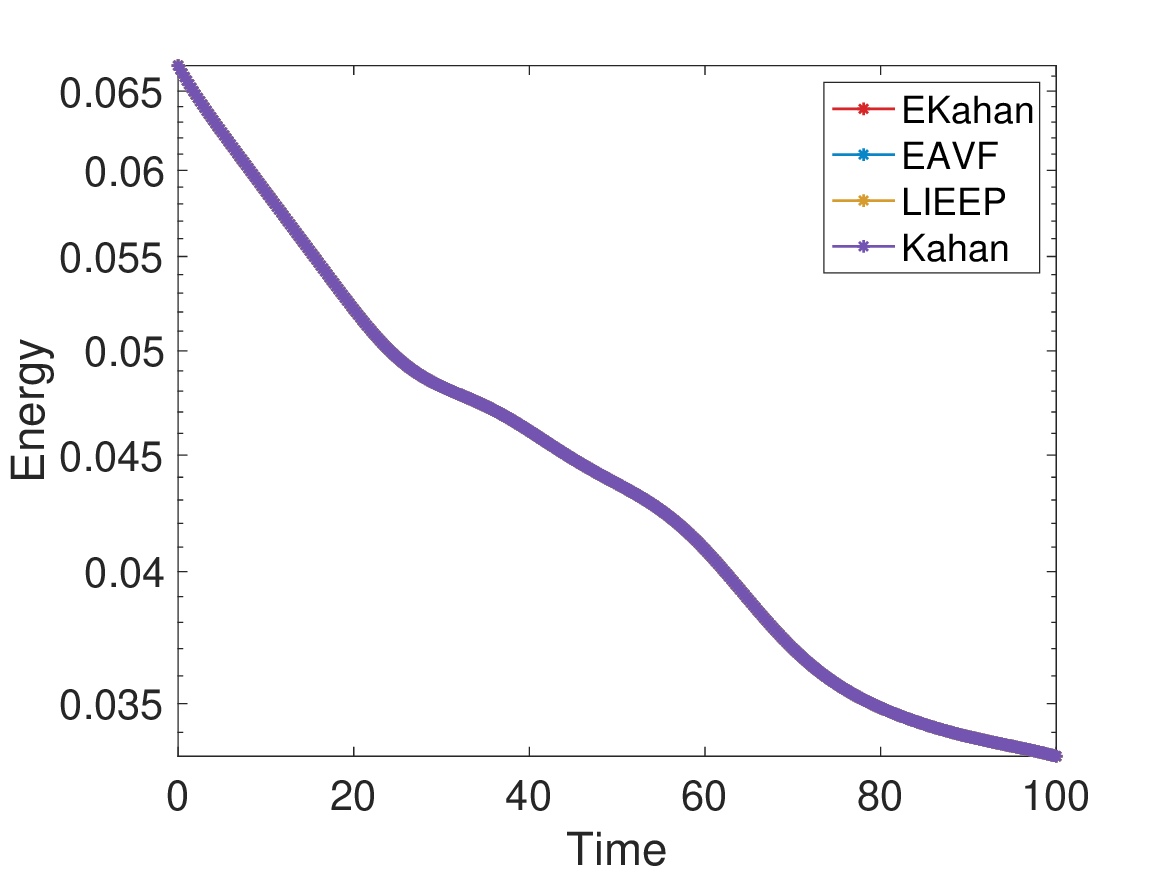}
		\caption[a]{$\gamma=0$, $\beta=2$}
		\label{fig:FPU energy_error_gamma=0_beta=2}
	\end{subfigure}
	\caption{FPU system with $p=1$, $T=100$, and $h_{2}=0.25$. Figure \ref{fig:FPU energy_error_gamma=0_beta=0}: the energy error defined in \eqref{Eerror-plot} computed by all schemes for the conservative system; Figure \ref{fig:FPU est}: the residual of equation \eqref{step-wise-E}  in Theorem \ref{theorem-energy} for the conservative system; Figure \ref{fig:FPU energy_error_gamma=0.1_beta=0} and Figure \ref{fig:FPU energy_error_gamma=0_beta=2}: the energy defined in Theorem \ref{theorem-energy} of all schemes for the two dissipative systems.}
	\label{fig:FPU energy_error}
\end{figure}

We consider three different settings of damping coefficients: 1) a conservative system with $\gamma = 0$ and $\beta = 0$, 2) a dissipative system with $\gamma = 0.1$ and $\beta = 0$, and 3) a dissipative system with $\gamma = 0$ and $\beta = 2$. We emphasize that our theoretical results are derived for the conservative Hamiltonian case, where the stepwise energy-increment identity in Theorem~\ref{theorem-energy} applies. The two dissipative settings are included as additional numerical tests to assess whether EKahan captures the expected energy-decay behavior when damping is present. Figure~\ref{fig:FPU energy_error} summarizes the observed energy behavior:
\begin{itemize}
	\item the  energy error of all methods stays bounded and oscillated  in the undamped case in Figure \ref{fig:FPU energy_error_gamma=0_beta=0},
	\item the per-step energy deviation under EKahan discretization shown in  Theorem \ref{theorem-energy} is validated in Figure  \ref{fig:FPU est} ,
	\item in the damped cases ($\gamma>0$ or $\beta>0$), the computed energy decreases over time, consistent with the dissipative nature of the underlying dynamics; see Figures~\ref{fig:FPU energy_error_gamma=0.1_beta=0} and \ref{fig:FPU energy_error_gamma=0_beta=2}.
\end{itemize}

\begin{figure}[H]
	\centering
	\begin{subfigure}[b]{0.4\textwidth}
		\includegraphics[width=\textwidth]{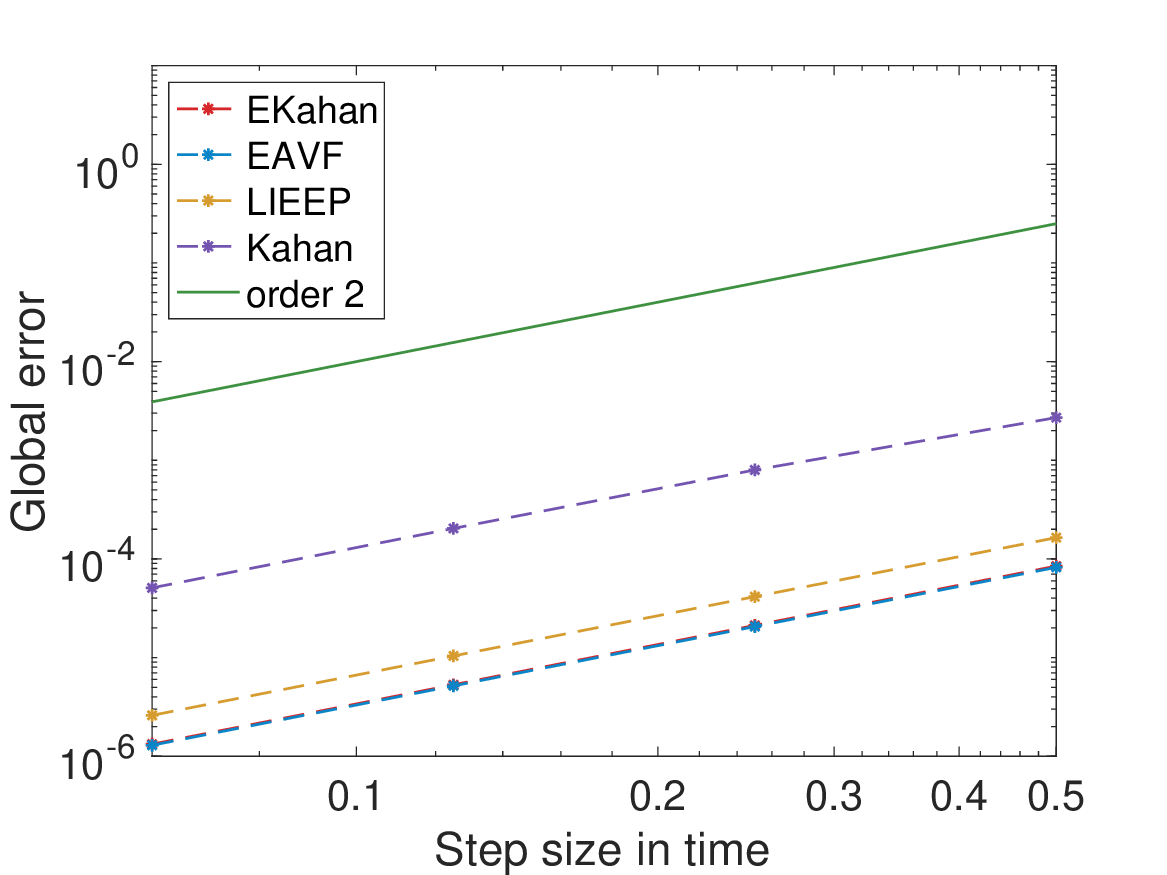}
		\caption[a]{order plot}
		\label{fig:FPU-p1-order}
	\end{subfigure}
	\begin{subfigure}[b]{0.4\textwidth}
		\includegraphics[width=\textwidth]{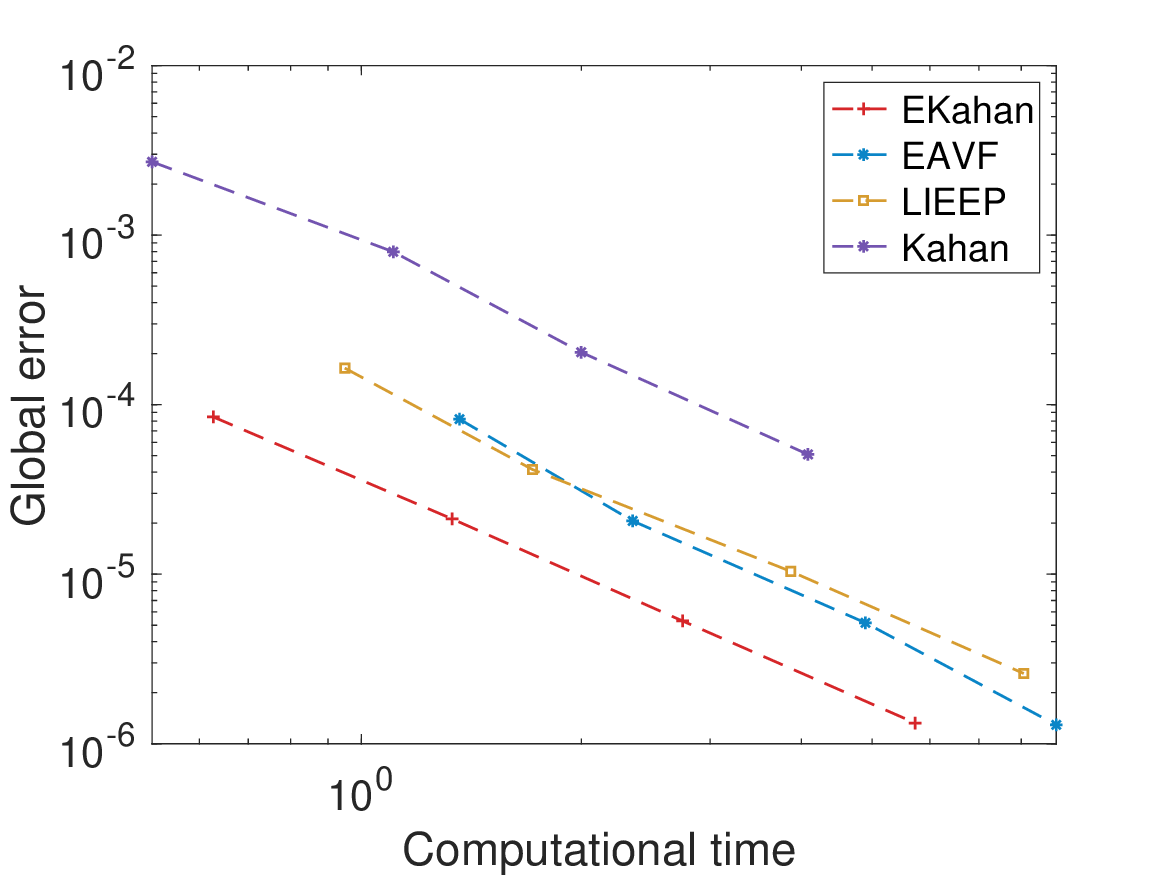}
		\caption[b]{efficiency plot}
		\label{fig:FPU-p1 efficiency}
	\end{subfigure}   
	\caption{Global error and computational time of different methods for the FPU system with $p=1$, $\gamma=0$, $\beta=0$, $T=100$ and different time step sizes $h_{i}$. }
	\label{fig:FPU-p1 order efficiency}
\end{figure}

\begin{figure}[H]
	\centering
	\begin{subfigure}[b]{0.4\textwidth}
		\includegraphics[width=\textwidth]{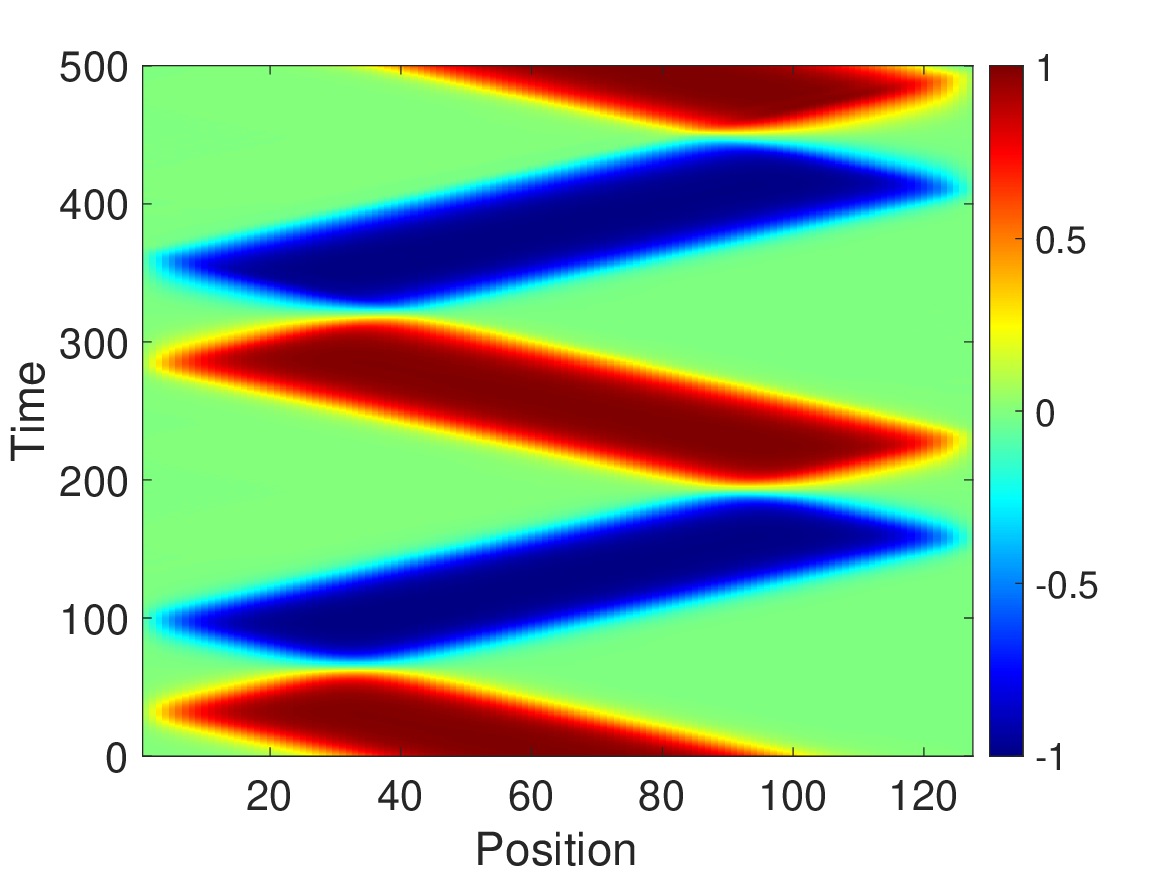}
		\caption[a]{$\gamma=0,\beta=0$}
		\label{fig:FPU solution gamma=beta=0}
	\end{subfigure}
	\begin{subfigure}[b]{0.4\textwidth}
		\includegraphics[width=\textwidth]{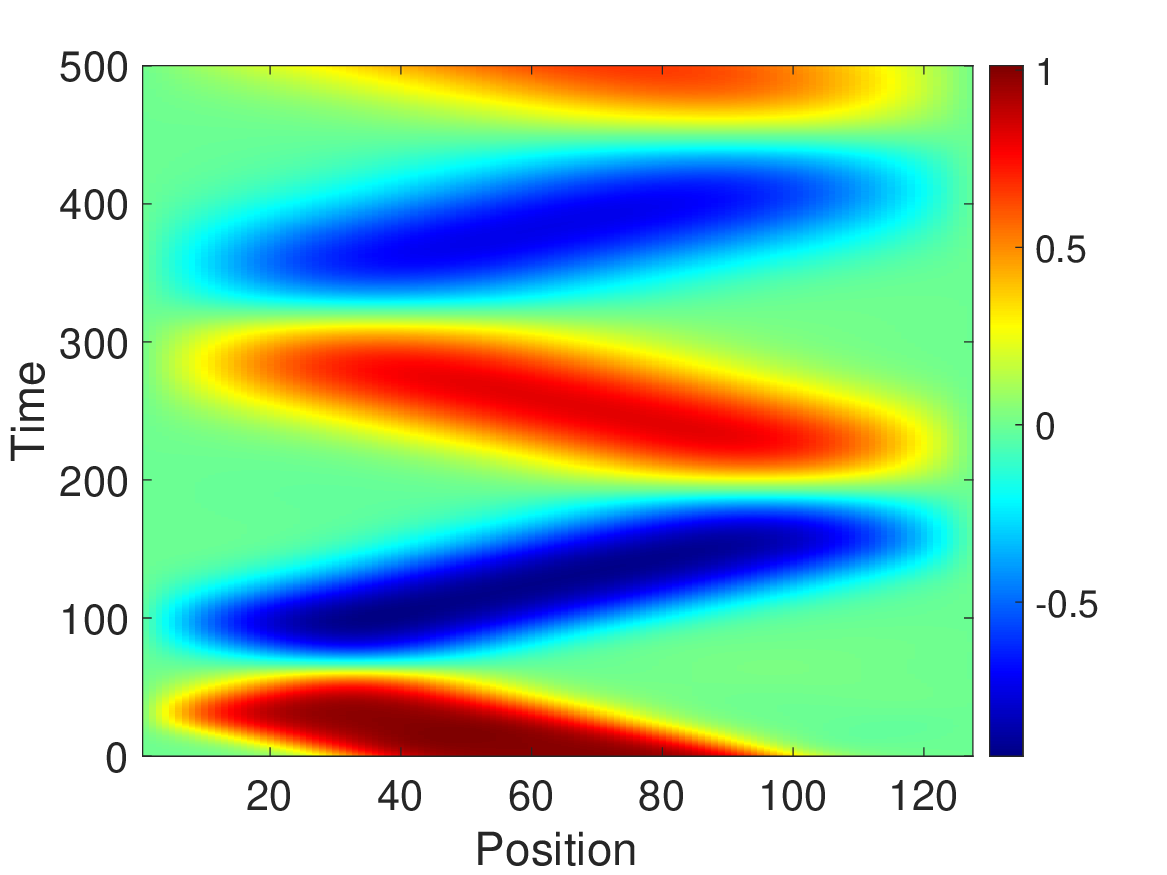}
		\caption[a]{$\gamma=0,\beta=2$}
		\label{fig:FPU solution gamma=0 beta=2}
	\end{subfigure}   
	\caption{ Numerical solutions of the FPU system computed over $t \in [0, 500]$ with time step $h_2 = 0.25$.}
	\label{fig:FPU solution}
\end{figure}

\begin{figure}[H]
	\centering
	\begin{subfigure}[b]{0.4\textwidth}
		\includegraphics[width=\textwidth]{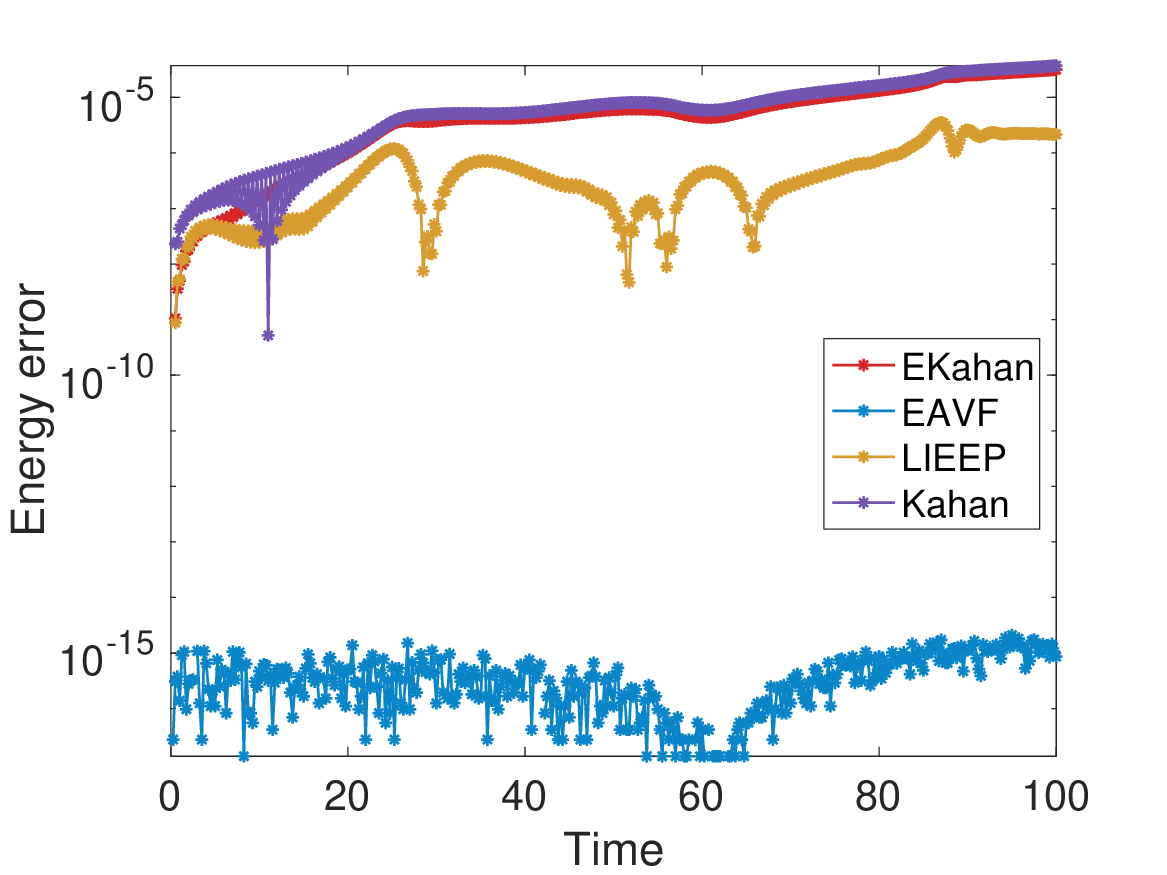}
		\caption[a]{}
		\label{fig:FPU p=2 energy error}
	\end{subfigure}
	\begin{subfigure}[b]{0.4\textwidth}
		\includegraphics[width=\textwidth]{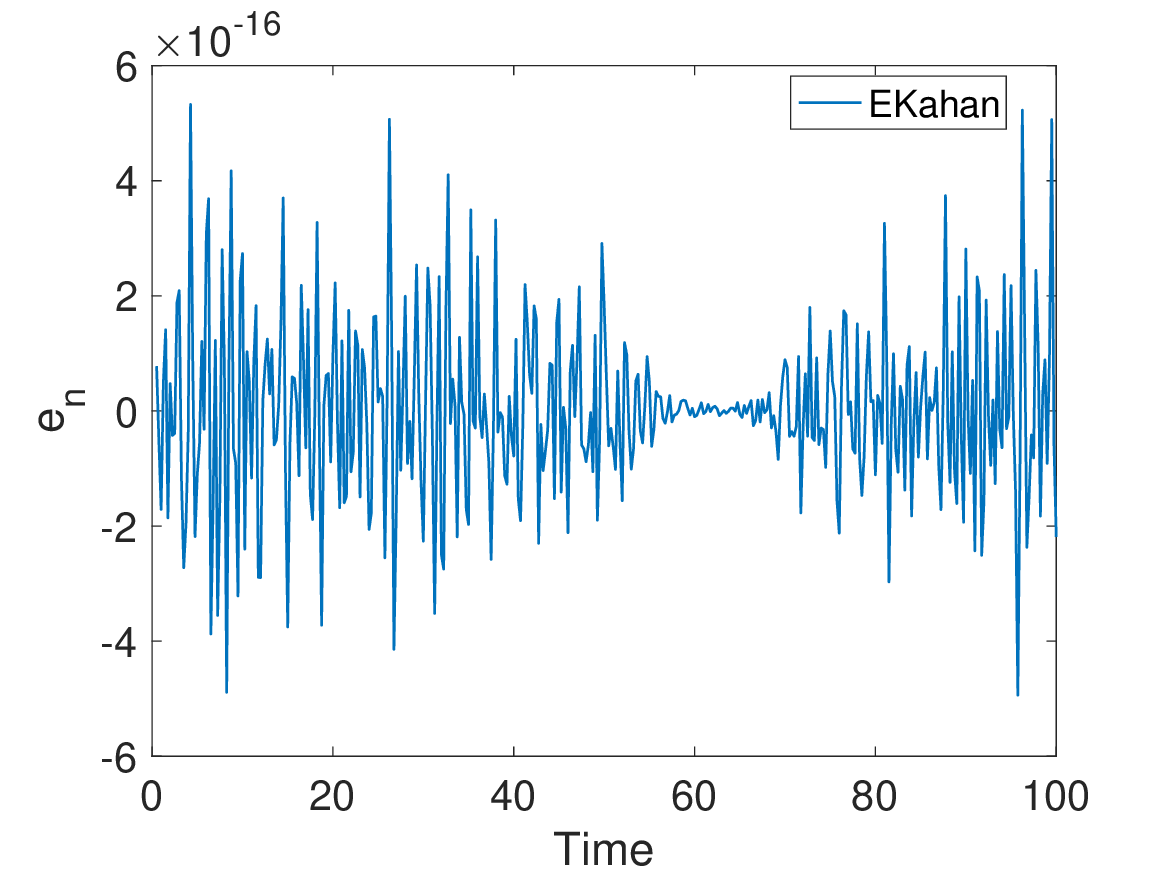}
		\caption[a]{}
		\label{fig:FPU p=2 est}
	\end{subfigure}   
	\caption{FPU system with $p=2$, $\beta=\gamma=0$, $T=100$ and $h_{2}=0.25$. Left: the energy error \eqref{Eerror-plot}  of all schemes; right: the residual of equation \eqref{higher-degree-EI} in Corollary \ref{corollary-EKahan-energy-high-order} by EKahan.}
	\label{fig:FPU p=2 energy}
\end{figure}

\begin{figure}[H]
	\centering
	\begin{subfigure}[b]{0.4\textwidth}
		\includegraphics[width=\textwidth]{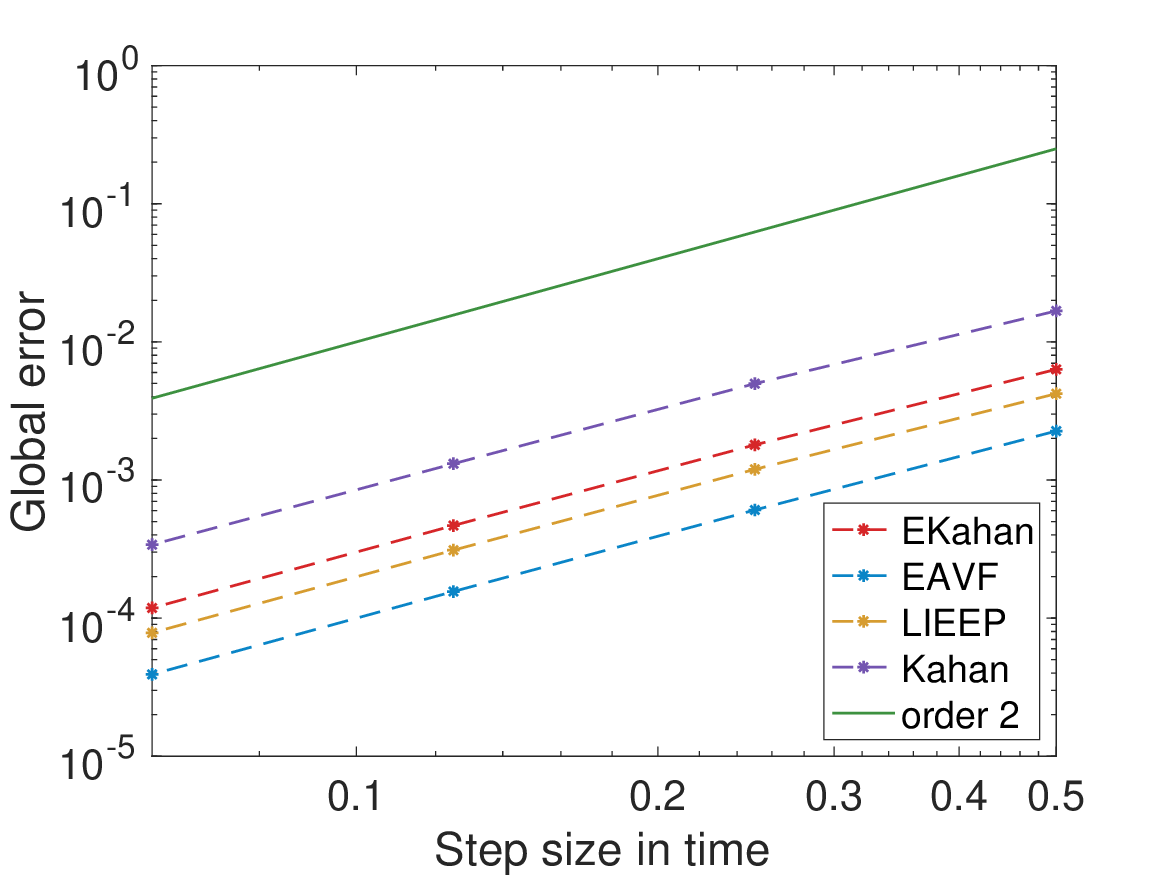}
		\caption[a]{order plot}
		\label{fig:FPU p=2 order}
	\end{subfigure}
	\begin{subfigure}[b]{0.4\textwidth}
		\includegraphics[width=\textwidth]{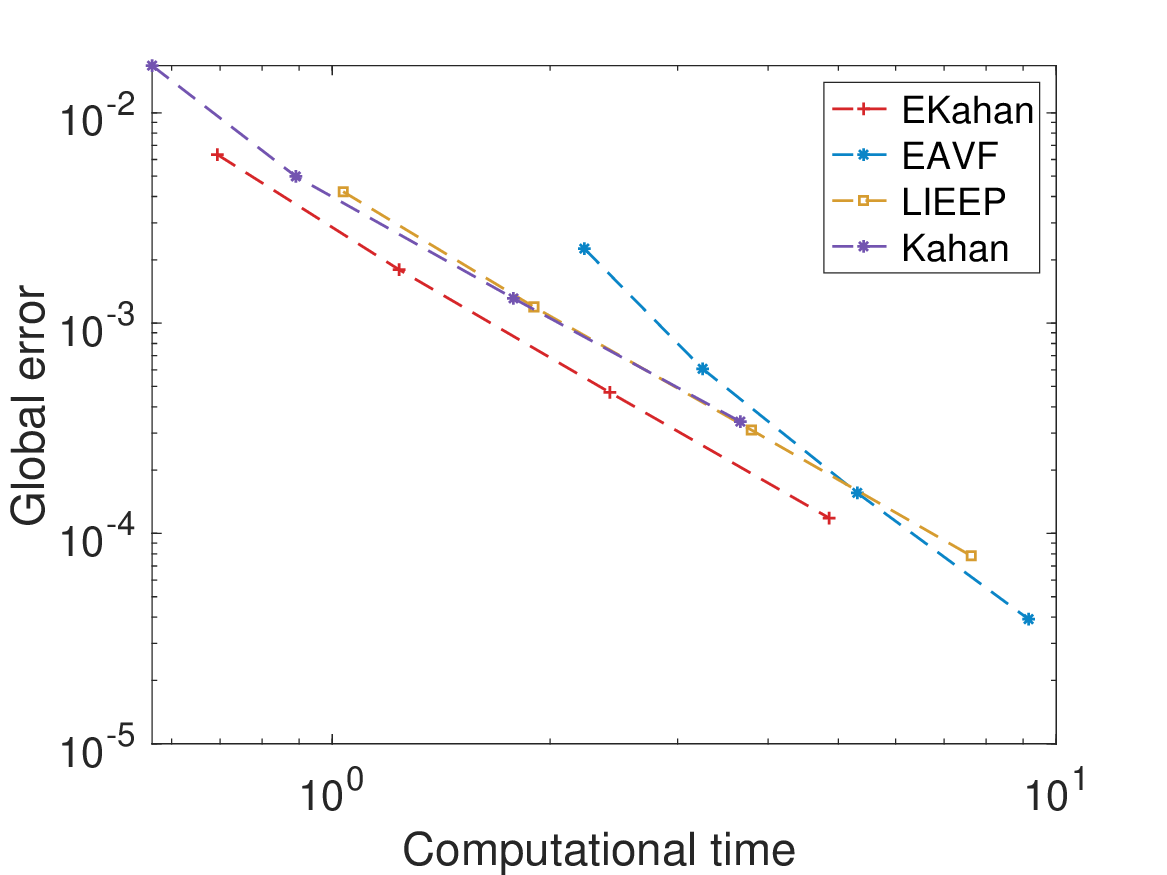}
		\caption[a]{efficiency plot}
		\label{fig:FPU p=2 efficiency}
	\end{subfigure}   
	\caption{Global error and computational time of different methods for the FPU system with $T=100$ and different time step sizes $h_{i}=1/2^{i}$, for $i=1,\dots,4$.}
	\label{fig:FPU p=2 order and efficiency}
\end{figure}

We test the convergence order and computational efficiency for all methods with those three settings of damping coefficients and the results are consistent. We report the results for the conservative equations in Figure \ref{fig:FPU-p1 order efficiency}.  Figure  \ref{fig:FPU-p1-order} 
shows that EKahan  is second-order and and Figure \ref{fig:FPU-p1 efficiency} shows that EKahan is more efficient than all other methods. In Figure \ref{fig:FPU solution}, we illustrate  the numerical
solutions given by EKahan method. Figure         \ref{fig:FPU solution gamma=beta=0} shows that EKahan captures the conservative property of the solution, exhibiting coherent traveling wave structures that preserve their shape and amplitude over time. Figure  \ref{fig:FPU solution gamma=0 beta=2} shows that EKahan  captures the decay of the wave profiles as energy dissipates internally. 
Figure \ref{fig:FPU energy_error_gamma=0.1_beta=0} and \ref{fig:FPU solution gamma=0 beta=2} demonstrate that although EKahan is designed for conservative equations, it can also capture the dissipative nature for the dissipative systems.  
\paragraph{\textbf{FPU system with a fourth-order Hamiltonian}} 
When \( p = 2 \), the potential function is a quartic homogeneous polynomial with the form
\[
U(y) = \sum_{j=0}^{N-1} \frac{\epsilon}{12} \left( \frac{u^{j+1} - u^j}{\Delta x} \right)^4.
\]
Define the spatial discrete difference by
\[
w^j = \frac{u^{j+1} - u^j}{\Delta x}, \quad j = 0, \ldots, N-1.
\]
Following equation \eqref{function-polarization}, the polarized  gradient  is given by
\begin{align*}
	\left(\nabla_{\mathrm{K}} U (y_n, y_{n+1}, y_{n+2})\right)^j 
	= \frac{\epsilon}{3 \Delta x} \left( w^{j-1}_n w^{j-1}_{n+1} w^{j-1}_{n+2} - w^j_n w^j_{n+1} w^j_{n+2} \right) 
\end{align*}
for \( j = 0, \ldots, N-1 \). The EKahan scheme then takes the form
\begin{equation*}
	y_{n+2} = e^{2h Q M} y_n + 2h \phi(2h Q M) Q \nabla_{\mathrm{K}} U (y_n, y_{n+1}, y_{n+2}).
	\label{eq:EKahan_update}
\end{equation*}
We employ identical initial boundary condition and computational parameters to the $p=1$ case, with the exception of setting $\epsilon = 100$. Figure \ref{fig:FPU p=2 energy error} illustrates that the
EAVF method preserves the discrete energy exactly and the energy errors of all linearly implicit methods remain bounded and oscillated. Figure \ref{fig:FPU p=2 est} verifies the formula of the energy deviation by EKahan method shown in Corollary \ref{corollary-EKahan-energy-high-order}.
In Figure \ref{fig:FPU p=2 order and efficiency}, we consider the global error  and the computational cost using different time step sizes. Figure \ref{fig:FPU p=2 order} demonstrates that EKahan method is second-order, and Figure \ref{fig:FPU p=2 efficiency} reveals its superior computational efficiency compared to all the other methods.

\paragraph{\textbf{Zakharov--Kuznetsov equation}}

The Zakharov--Kuznetsov (ZK) equation is a two-dimensional generalization of the classical Korteweg--de Vries (KdV) equation~\cite{ZK1974}, modeling the propagation of nonlinear waves in magnetized plasmas and stratified fluids. In this study, we consider the specific form
\begin{equation*}
	u_t + u^p u_x + u_{xxx} + u_{xyy} = 0,
\end{equation*}
for $(x, y) \in [0, L] \times [0, L]$ and $t \in [0, T]$. Among its conserved quantities, we consider the following energy functional
\begin{equation*}
	\mathcal{H}(t) = \frac{1}{2} \int_{\mathbb{R}^2} (u_x^2 + u_y^2) \, \mathrm{d}x \, \mathrm{d}y - \frac{1}{(p+1)(p+2)} \int_{\mathbb{R}^2} u^{p+2} \, \mathrm{d}x \, \mathrm{d}y.
\end{equation*}
For $p = 1$, the equation admits the following Hamiltonian structure:
\begin{equation}\label{ZK-pde}
	u_t = -\partial_x \frac{\delta \mathcal{H}}{\delta u},\quad \frac{\delta \mathcal{H}}{\delta u} = \frac{1}{2} u^2 + u_{xx} + u_{yy}.
\end{equation}

To discretize the problem, we use the following discrete energy 
\begin{equation*}
	H(u) = \sum_{j=0}^{N_y} \sum_{i=0}^{N_x} \left( \frac{(\delta_x^+ u^{i,j})^2 + (\delta_x^- u^{i,j})^2}{4} + \frac{(\delta_y^+ u^{i,j})^2 + (\delta_y^- u^{i,j})^2}{4}  + \frac{(u^{i,j})^3}{6} \right) \Delta x \Delta y.
\end{equation*}
Here, $\delta_x^+$ and $\delta_x^-$, similarly $\delta_y^+$ and $\delta_y^-$, denote the forward and backward finite difference operators in the $x$- and $y$-directions, respectively. First-order and second-order central difference operators are employed to approximate the first and second spatial derivatives in \eqref{ZK-pde}, respectively, and the corresponding  matrices are denoted by $D_x^{(1)}, D_x^{(2)}, D_y^{(1)},$ and $D_y^{(2)}$. Using these notations, the semi-discrete system can be compactly written as
\begin{equation*}
	\dot{U} = -D_x^{(1)} \left( \frac{1}{2} U^2 + D_x^{(2)} U + D_y^{(2)} U \right),
\end{equation*}
where the vector $U \in \mathbb{R}^{(N_x-1)\times(N_y-1)}$  is flattened from  $u_{i,j}$ by stacking columns successively, giving priority to the $y$-direction, and $U^2$ denotes the element-wise square of $U$. Here, $N_x$ and $N_y$ are the grid numbers in each spatial direction.

Numerical simulations are carried out on the spatial domain $[0, L] \times [0, L]$ with $L = 6$ and the temporal domain $[0, T]$ with $T = 8$, subject to periodic boundary conditions. The spatial discretization employs a uniform grid with $N = 33$ points in each spatial direction, resulting in a mesh size of $\Delta x = \Delta y = L / (N-1)$. We consider different time steps $h_i = 0.01/2^{i+1}, \; i = 1,\cdots, 4$, and the initial condition
\begin{equation*}
	u(0,x,y)=\sqrt{2}[\sin(\frac{2\pi x}{L})+\frac{1}{\sqrt{2}} \cos (\frac{4\pi x}{L}+\frac{\pi}{4})]\times[\cos(\frac{2\pi y}{L})+\frac{1}{\sqrt{2}} \cos (\frac{4\pi y}{L}+\frac{\pi}{3})].
\end{equation*}

\begin{figure}[H]
	\centering
	\begin{subfigure}[b]{0.4\textwidth}
		\includegraphics[width=\textwidth]{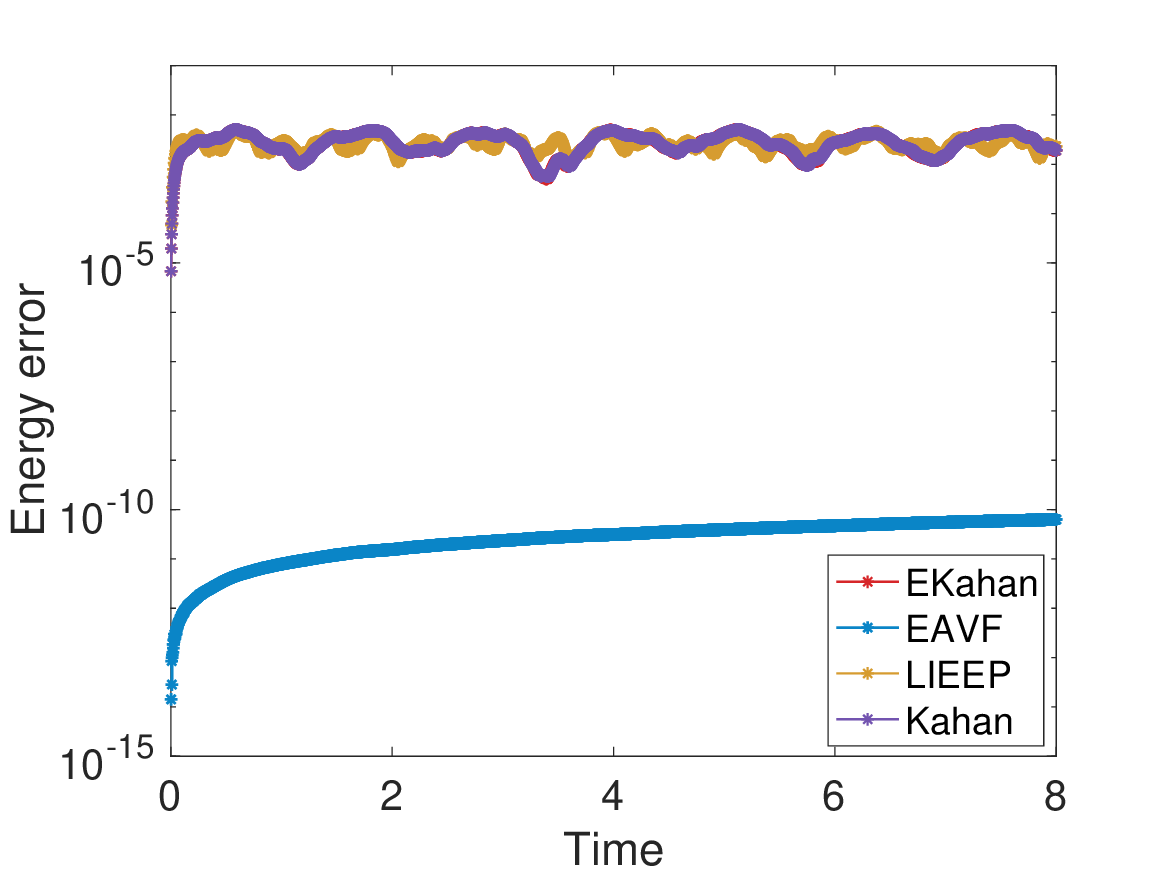}
		\caption[a]{}
		\label{fig:ZK energy error}
	\end{subfigure}
	\begin{subfigure}[b]{0.4\textwidth}
		\includegraphics[width=\textwidth]{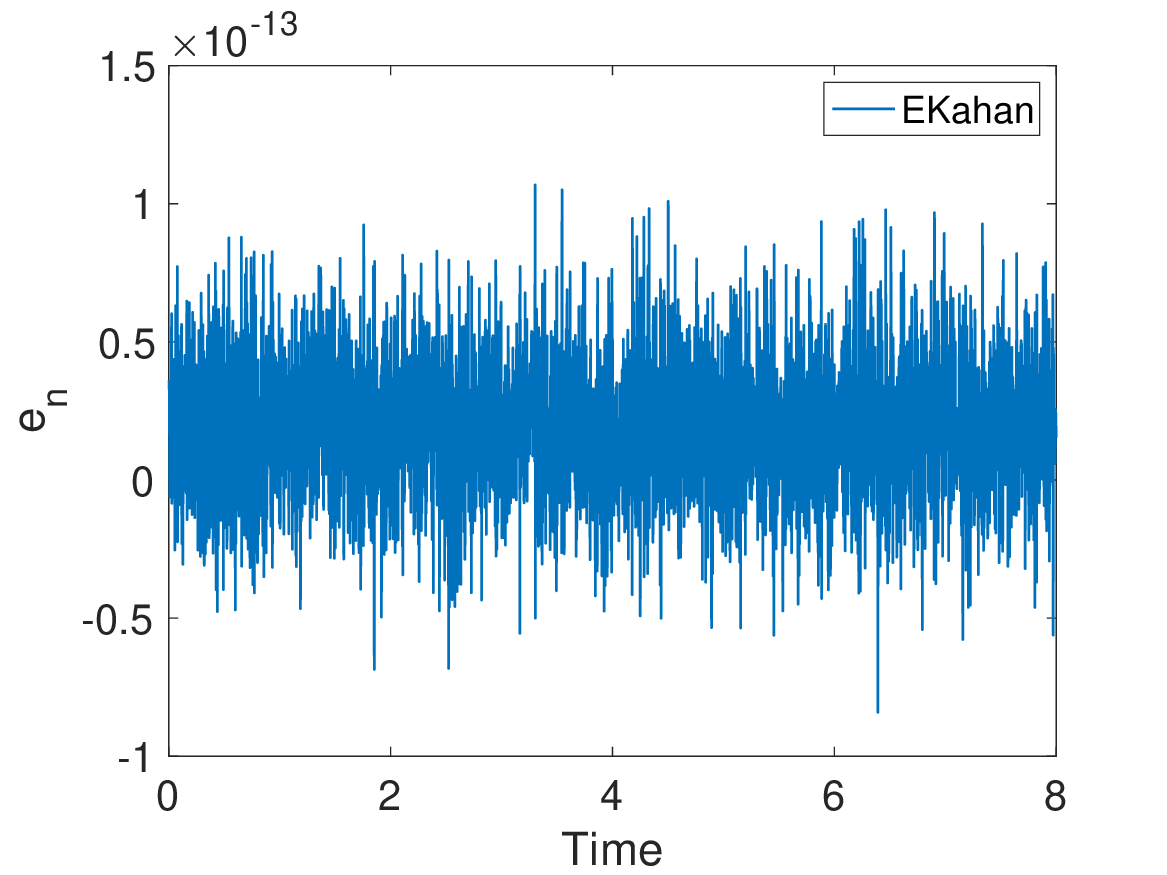}
		\caption[a]{}
		\label{fig:ZK energy est}
	\end{subfigure}   
	\caption{Zakharov-Kuznetsov equation with $T=8$ and $h_{1}=0.0025$.
		Left: the energy error defined in \eqref{Eerror-plot}; right: the residual of equation \eqref{step-wise-E}  in Theorem \ref{theorem-energy}.}
	\label{fig:ZK energy error and energy est}
\end{figure}

Figure \ref{fig:ZK energy error} demonstrates that  the energy errors of all linearly implicit methods remain bounded and oscillated. Figure \ref{fig:ZK energy est} confirms the stepwise energy-increment identity for EKahan established in Theorem~\ref{theorem-energy}.
\begin{figure}[H]
	\centering
	\begin{subfigure}[b]{0.4\textwidth}
		\includegraphics[width=\textwidth]{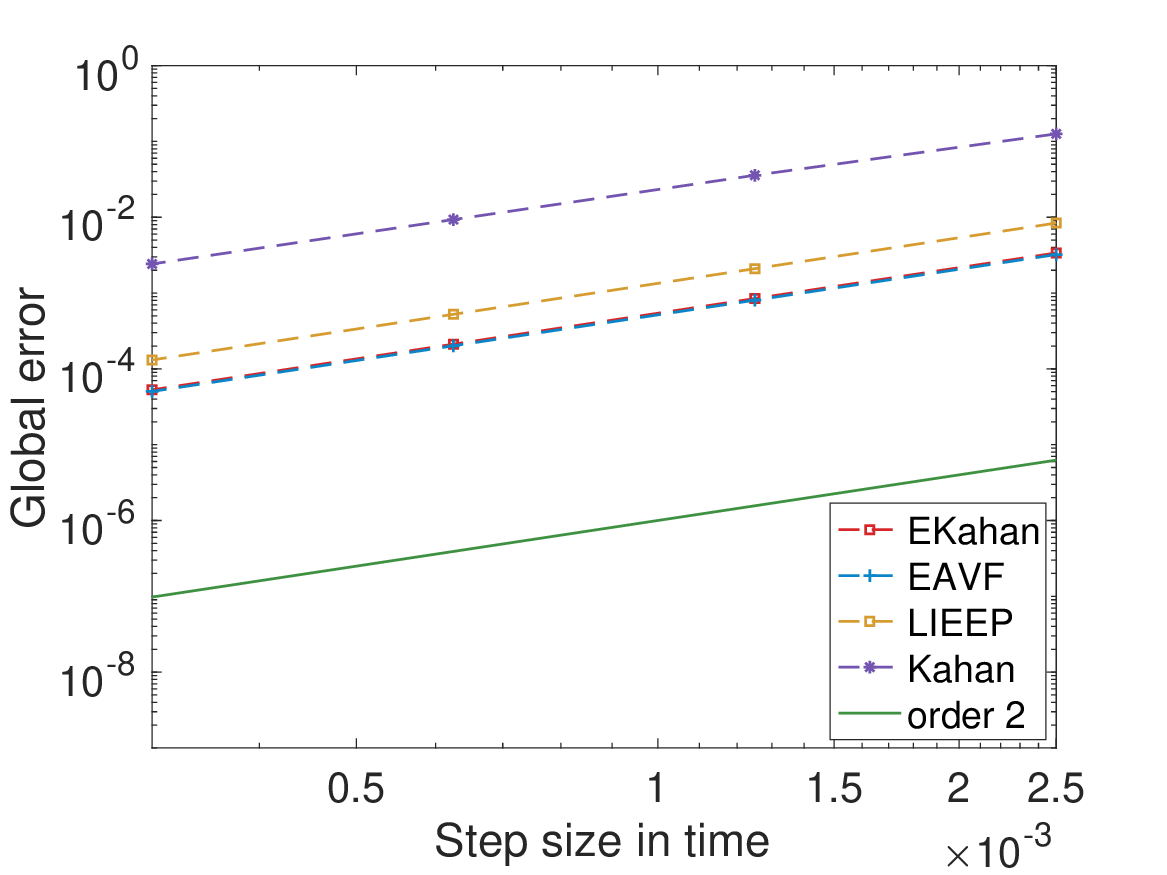}
		\caption[a]{order plot}
		\label{fig:ZK order}
	\end{subfigure}
	\begin{subfigure}[b]{0.4\textwidth}
		\includegraphics[width=\textwidth]{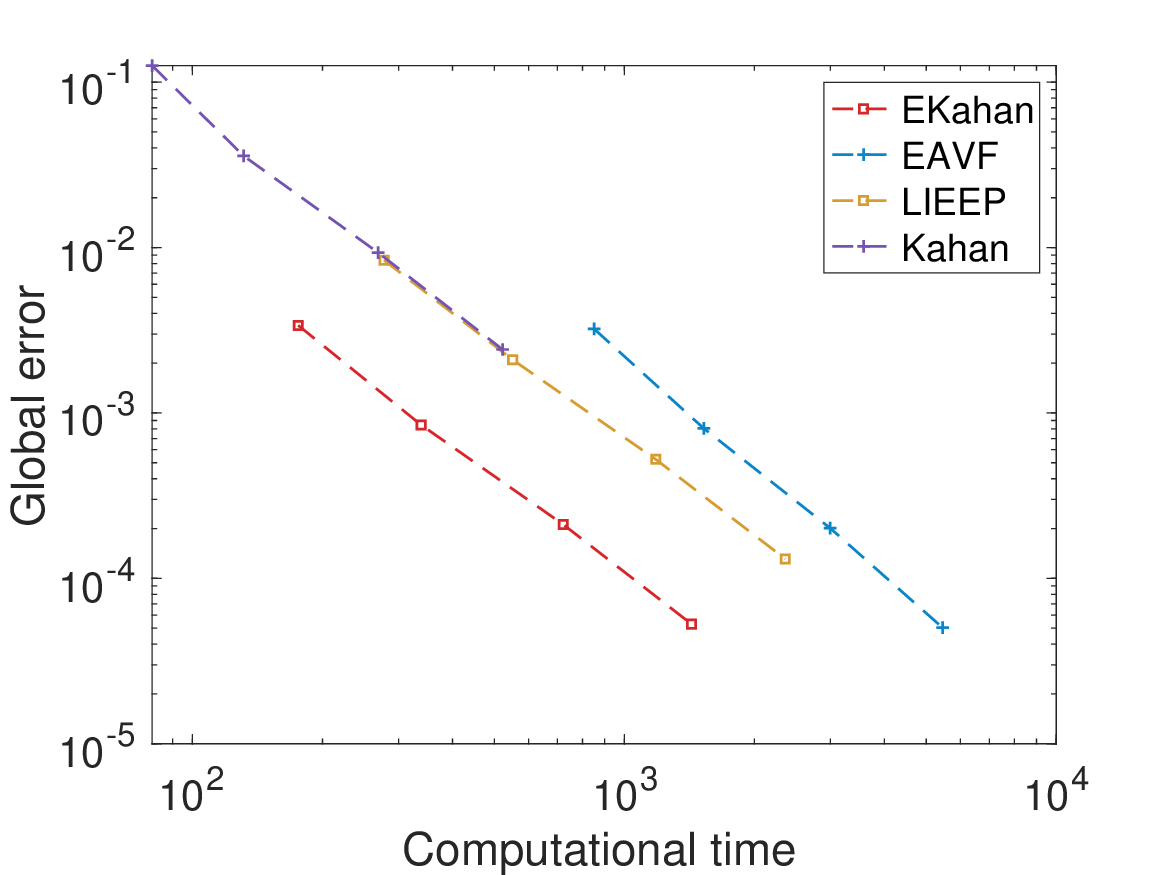}
		\caption[a]{efficiency plot}
		\label{fig:ZK efficiency}
	\end{subfigure}   
	\caption{Global error and computational time of different methods for the Zakharov-Kuznetsov equation with $T=8$ and different time step sizes $h_{i}$.}
	\label{fig:ZK order and efficiency}
\end{figure}
\begin{figure}[H]
	\centering
	\begin{subfigure}[b]{0.4\textwidth}
		\includegraphics[width=\textwidth]{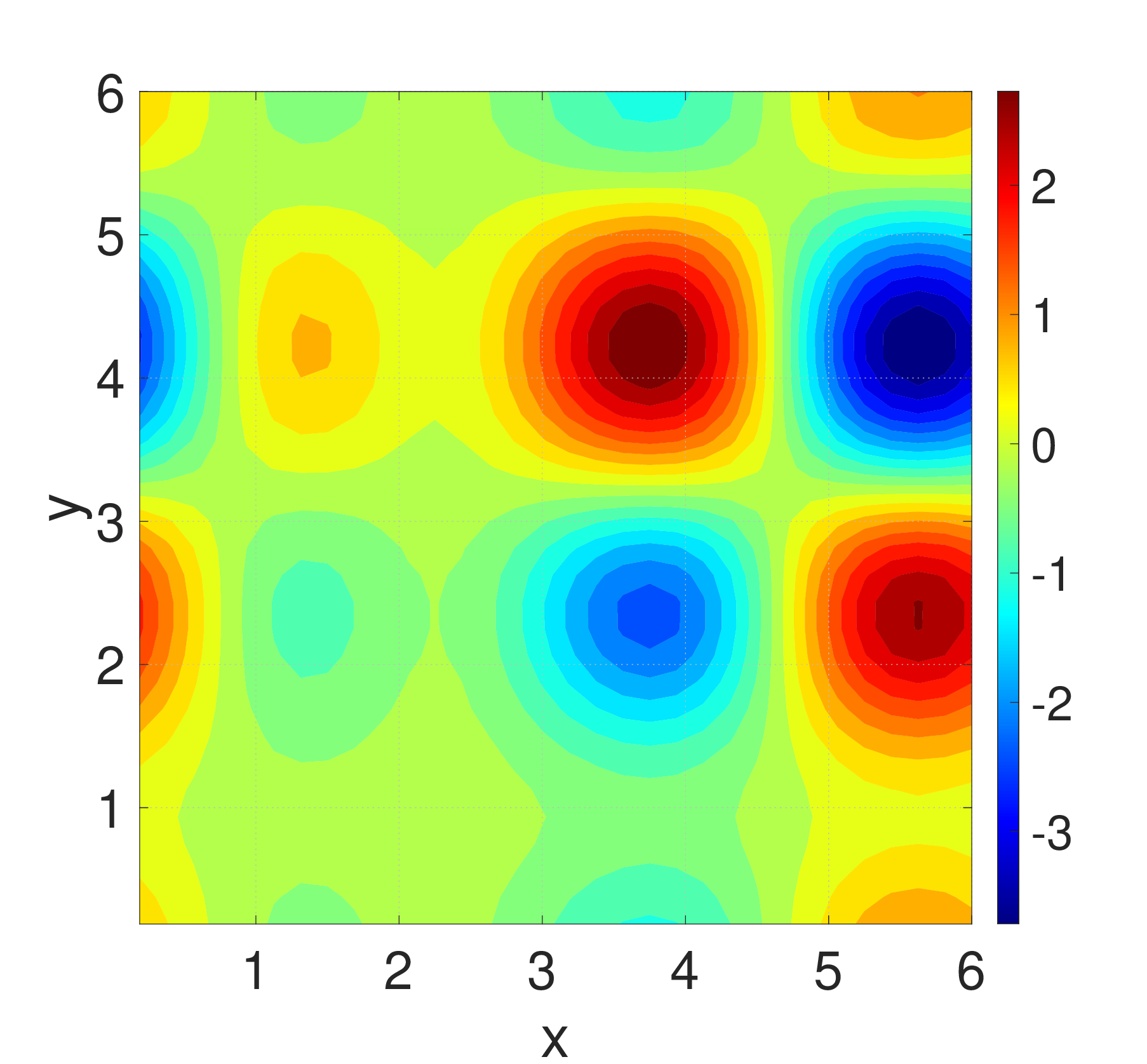}
		\caption[t=0]{}
		\label{fig:ZK solution initial}
	\end{subfigure}
	\begin{subfigure}[b]{0.4\textwidth}
		\includegraphics[width=\textwidth]{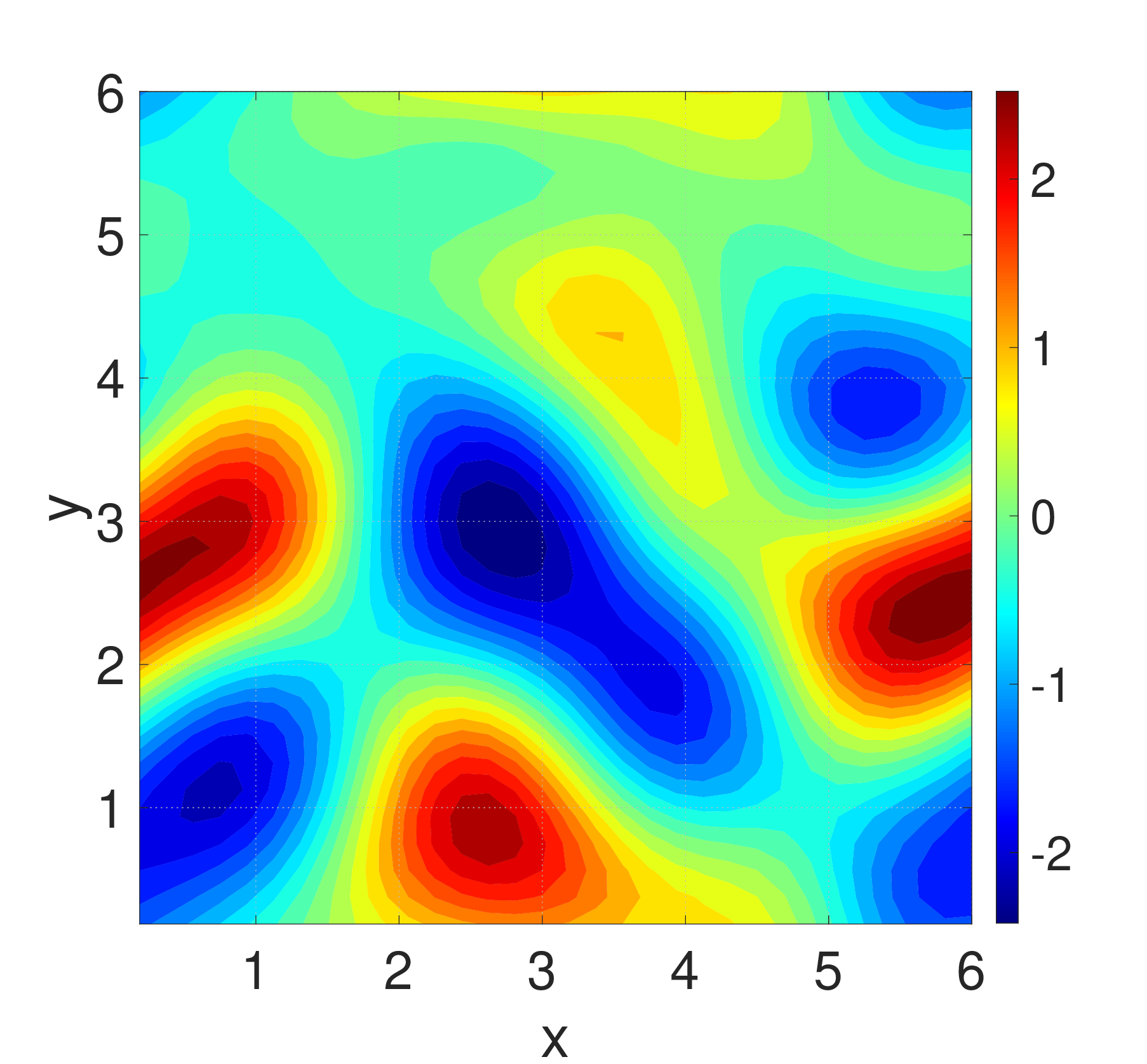}
		\caption[t=8]{}
		\label{fig:ZK solution end}
	\end{subfigure}   
	\caption{Evolution to bell-shaped solitons in Zakharov-Kuznetsov equation, resolved with time step $h_1 = 0.0025$. Left: initial condition; right: soliton structure at $t = 8$. }
	\label{fig:ZK solution}
\end{figure}
In Figure \ref{fig:ZK order and efficiency}, we consider the global error and the computational cost using different step sizes $h_{i}$. Figure \ref{fig:ZK order} demonstrates the second-order convergence
of the EKahan method, and Figure \ref{fig:ZK efficiency} reveals its superior computational efficiency compared to  the other methods.  Figure  \ref{fig:ZK solution} illustrates the solution of the Zakharov-Kuznetsov equation. In Figure \ref{fig:ZK solution initial}, we observe that  the initial condition exhibits a regular, periodic pattern composed of multiple interacting modes. As time progresses to 
$t=8$, the solution evolves into a more localized, bell-shaped structure, reflecting the nonlinear dispersive nature of the equation; as shown in Figure \ref{fig:ZK solution end}.  Despite the significant structural changes in the solution profile, the total energy of the system remains nearly constant throughout the simulation, demonstrating the excellent long-time stability.

\section{Conclusion}
In this work, we proposed a linearly implicit structure-preserving exponential integrator for semilinear Hamiltonian systems with polynomial nonlinearities by combining the Kahan-type discretization with exponential integration techniques. {\color{black}Although not all geometric properties of Kahan’s method are exactly inherited by EKahan, the resulting scheme balances computational efficiency with key geometric features, including symmetry and exact stepwise energy-increment identities. Together with our second-order convergence result, these identities provide a computable near conservation of energy characterization in terms of numerical increments, rather than implying exact energy conservation or symplecticity for general nonlinear systems. Compared to the fully implicit energy-preserving exponential integrator EAVF method, EKahan requires only the solution of a single linear system per time step, significantly reducing computational costs.}

Through extensive numerical experiments on benchmark systems, we verified the method’s long-term accuracy and computational competitiveness. In particular, EKahan was shown to outperform not only the fully implicit scheme but also other linearly implicit methods, LIEEP and Kahan method, highlighting its enhanced efficiency among structure-preserving approaches. {\color{black}We also tested the method in representative damped regimes to assess its robustness beyond the conservative setting, where EKahan captures the expected energy-decay trend observed in the dynamics.} Furthermore, we provided  the generalization of EKahan to problems involving higher-degree polynomial Hamiltonians, thereby broadening its applicability to a wider class of physical models.

Overall, the EKahan method offers a compelling tool for the long-time integration of semilinear Hamiltonian systems. Future work may explore extending this framework to adaptive step-size strategies, and constructing higher-order efficient structure-preserving exponential integrators.


\begin{thebibliography}{10}
	\providecommand{\url}[1]{{#1}}
	\providecommand{\urlprefix}{URL }
	\expandafter\ifx\csname urlstyle\endcsname\relax
	\providecommand{\doi}[1]{DOI~\discretionary{}{}{}#1}\else
	\providecommand{\doi}{DOI~\discretionary{}{}{}\begingroup
		\urlstyle{rm}\Url}\fi
	
	\bibitem{ambati2015review}
	Ambati, M., Gerasimov, T., De~Lorenzis, L.: A review on phase-field models of
	brittle fracture and a new fast hybrid formulation.
	\newblock Comput. Mech. \textbf{55}(2), 383--405 (2015).
	\newblock \doi{10.1007/s00466-014-1109-y}
	
	\bibitem{1206044}
	Berland, H., Skaflestad, B., Wright, W.M.: {EXPINT}---{A MATLAB} package for
	exponential integrators.
	\newblock ACM Trans. Math. Softw. \textbf{33}(1), 4–es (2007).
	\newblock \doi{10.1145/1206040.1206044}
	
	\bibitem{Cary}
	Cary, J.R., Brizard, A.J.: Hamiltonian theory of guiding-center motion.
	\newblock Rev. Modern Phys. \textbf{81}(2), 693--738 (2009).
	\newblock \doi{10.1103/RevModPhys.81.693}
	
	\bibitem{Celledoni2008}
	Celledoni, E., Cohen, D., Owren, B.: Symmetric exponential integrators with an
	application to the cubic {S}chr\"odinger equation.
	\newblock Found. Comput. Math. \textbf{8}(3), 303--317 (2008).
	\newblock \doi{10.1007/s10208-007-9016-7}
	
	\bibitem{Celledoni2012}
	Celledoni, E., Grimm, V., McLachlan, R.I., McLaren, D.I., O'Neale, D., Owren,
	B., Quispel, G.R.W.: Preserving energy resp. dissipation in numerical {PDE}s
	using the ``average vector field'' method.
	\newblock J. Comput. Phys. \textbf{231}(20), 6770--6789 (2012).
	\newblock \doi{10.1016/j.jcp.2012.06.022}
	
	\bibitem{MR3452144}
	Celledoni, E., McLachlan, R.I., McLaren, D.I., Owren, B., Quispel, G.R.W.:
	Discretization of polynomial vector fields by polarization.
	\newblock Proc. A. \textbf{471}(2184), 20150390, 10 (2015).
	\newblock \doi{10.1098/rspa.2015.0390}
	
	\bibitem{celledoni2012geometric}
	Celledoni, E., McLachlan, R.I., Owren, B., Quispel, G.R.W.: Geometric
	properties of {K}ahan's method.
	\newblock J. Phys. A \textbf{46}(2), 025201, 12 (2013).
	\newblock \doi{10.1088/1751-8113/46/2/025201}
	
	\bibitem{DENG2022166}
	Deng, S., Li, J.: A uniformly accurate exponential wave integrator {F}ourier
	pseudo-spectral method with energy-preservation for long-time dynamics of the
	nonlinear {K}lein-{G}ordon equation.
	\newblock Appl. Numer. Math. \textbf{178}, 166--191 (2022).
	\newblock \doi{10.1016/j.apnum.2022.03.019}
	
	\bibitem{eidnes2019linearly}
	Eidnes, S., Li, L., Sato, S.: Linearly implicit structure-preserving schemes
	for {H}amiltonian systems.
	\newblock J. Comput. Appl. Math. \textbf{387}, Paper No. 112489, 12 (2021).
	\newblock \doi{10.1016/j.cam.2019.112489}
	
	\bibitem{4376203}
	Fermi, E., Pasta, P., Ulam, S., Tsingou, M.: {S}tudies of the nonlinear
	problem.
	\newblock Tech. rep., Los Alamos National Laboratory (LANL), Los Alamos, NM
	(United States) (1955).
	\newblock \doi{10.2172/4376203}
	
	\bibitem{MR1128990}
	Fordy, A.P.: The {H}\'enon-{H}eiles system revisited.
	\newblock Phys. D \textbf{52}(2-3), 204--210 (1991).
	\newblock \doi{10.1016/0167-2789(91)90122-P}
	
	\bibitem{fu2022high}
	Fu, Y., Hu, D., Xu, Z.: High-order explicit conservative exponential integrator
	schemes for fractional {H}amiltonian {PDE}s.
	\newblock Appl. Numer. Math. \textbf{172}, 315--331 (2022).
	\newblock \doi{10.1016/j.apnum.2021.10.011}
	
	\bibitem{fu2022arbitrary}
	Fu, Y., Hu, D., Zhang, G.: Arbitrary high-order exponential integrators
	conservative schemes for the nonlinear {G}ross-{P}itaevskii equation.
	\newblock Comput. Math. Appl. \textbf{121}, 102--114 (2022).
	\newblock \doi{10.1016/j.camwa.2022.07.004}
	
	\bibitem{hairer2010energy}
	Hairer, E.: Energy-preserving variant of collocation methods.
	\newblock JNAIAM. J. Numer. Anal. Ind. Appl. Math. \textbf{5}(1-2), 73--84
	(2010)
	
	\bibitem{hairer2006geometric}
	Hairer, E., Lubich, C., Wanner, G.: Geometric numerical integration:
	structure-preserving algorithms for ordinary differential equations, vol.~31.
	\newblock Springer Science \& Business Media (2006)
	
	\bibitem{MR2652783}
	Hochbruck, M., Ostermann, A.: Exponential integrators.
	\newblock Acta Numer. \textbf{19}, 209--286 (2010).
	\newblock \doi{10.1017/S0962492910000048}
	
	\bibitem{jiang2022high}
	Jiang, C., Cui, J., Qian, X., Song, S.: High-order linearly implicit
	structure-preserving exponential integrators for the nonlinear
	{S}chr\"odinger equation.
	\newblock J. Sci. Comput. \textbf{90}(1), Paper No. 66, 27 (2022).
	\newblock \doi{10.1007/s10915-021-01739-x}
	
	\bibitem{jiang2020linearly}
	Jiang, C., Wang, Y., Cai, W.: A linearly implicit energy-preserving exponential
	integrator for the nonlinear {K}lein-{G}ordon equation.
	\newblock J. Comput. Phys. \textbf{419}, 109690, 18 (2020).
	\newblock \doi{10.1016/j.jcp.2020.109690}
	
	\bibitem{Ju2022}
	Ju, L., Li, X., Qiao, Z.: Generalized {SAV}-exponential integrator schemes for
	{A}llen-{C}ahn type gradient flows.
	\newblock SIAM J. Numer. Anal. \textbf{60}(4), 1905--1931 (2022).
	\newblock \doi{10.1137/21M1446496}
	
	\bibitem{kato1987nonlinear}
	Kato, T.: On nonlinear {S}chr{\"o}dinger equations.
	\newblock In: Annales de l'I.H.P. Physique th{\'e}orique, vol.~46, pp. 113--129
	(1987)
	
	\bibitem{li2022energy}
	Li, J.: Energy-preserving exponential integrator {F}ourier pseudo-spectral
	schemes for the nonlinear {D}irac equation.
	\newblock Appl. Numer. Math. \textbf{172}, 1--26 (2022).
	\newblock \doi{10.1016/j.apnum.2021.09.006}
	
	\bibitem{li2022new}
	Li, L.: A new symmetric linearly implicit exponential integrator preserving
	polynomial invariants or {L}yapunov functions for conservative or dissipative
	systems.
	\newblock J. Comput. Phys. \textbf{449}, Paper No. 110800, 13 (2022).
	\newblock \doi{10.1016/j.jcp.2021.110800}
	
	\bibitem{Wuxinyuan2016}
	Li, Y.W., Wu, X.: Exponential integrators preserving first integrals or
	{L}yapunov functions for conservative or dissipative systems.
	\newblock SIAM J. Sci. Comput. \textbf{38}(3), A1876--A1895 (2016).
	\newblock \doi{10.1137/15M1023257}
	
	\bibitem{liu2020}
	Liu, Z., Li, X.: The exponential scalar auxiliary variable ({E}-{SAV}) approach
	for phase field models and its explicit computing.
	\newblock SIAM J. Sci. Comput. \textbf{42}(3), B630--B655 (2020).
	\newblock \doi{10.1137/19M1305914}
	
	\bibitem{matsuo2001dissipative}
	Matsuo, T., Furihata, D.: Dissipative or conservative finite-difference schemes
	for complex-valued nonlinear partial differential equations.
	\newblock J. Comput. Phys. \textbf{171}(2), 425--447 (2001).
	\newblock \doi{10.1006/jcph.2001.6775}
	
	\bibitem{mei2019exponential}
	Mei, L., Huang, L., Huang, S.: Exponential integrators with quadratic energy
	preservation for linear {P}oisson systems.
	\newblock J. Comput. Phys. \textbf{387}, 446--454 (2019).
	\newblock \doi{10.1016/j.jcp.2019.03.005}
	
	\bibitem{MEI2021110429}
	Mei, L., Huang, L., Wu, X.: Energy-preserving exponential integrators of
	arbitrarily high order for conservative or dissipative systems with highly
	oscillatory solutions.
	\newblock J. Comput. Phys. \textbf{442}, Paper No. 110429, 22 (2021).
	\newblock \doi{10.1016/j.jcp.2021.110429}
	
	\bibitem{mei2022energy}
	Mei, L., Huang, L., Wu, X.: Energy-preserving continuous-stage exponential
	{R}unge-{K}utta integrators for efficiently solving {H}amiltonian systems.
	\newblock SIAM J. Sci. Comput. \textbf{44}(3), A1092--A1115 (2022).
	\newblock \doi{10.1137/21M1412475}
	
	\bibitem{MEI2022110822}
	Mei, L., Huang, L., Wu, X.: A unified framework for the study of high-order
	energy-preserving integrators for solving {P}oisson systems.
	\newblock J. Comput. Phys. \textbf{450}, Paper No. 110822, 27 (2022).
	\newblock \doi{10.1016/j.jcp.2021.110822}
	
	\bibitem{sato2024high}
	Sato, S.: High-order linearly implicit exponential integrators conserving
	quadratic invariants with application to scalar auxiliary variable approach.
	\newblock Numerical Algorithms \textbf{96}(3), 1295--1329 (2024)
	
	\bibitem{shen2019new}
	Shen, J., Xu, J., Yang, J.: A new class of efficient and robust energy stable
	schemes for gradient flows.
	\newblock SIAM Rev. \textbf{61}(3), 474--506 (2019).
	\newblock \doi{10.1137/17M1150153}
	
	\bibitem{wang2019volume}
	Wang, B., Wu, X.: Volume-preserving exponential integrators and their
	applications.
	\newblock J. Comput. Phys. \textbf{396}, 867--887 (2019).
	\newblock \doi{10.1016/j.jcp.2019.07.026}
	
	\bibitem{Wuxinyuan2018}
	Wu, X., Wang, B.: Recent developments in structure-preserving algorithms for
	oscillatory differential equations.
	\newblock Science Press Beijing, Beijing; Springer, Singapore (2018).
	\newblock \doi{10.1007/978-981-10-9004-2}
	
	\bibitem{wu2012explicit}
	Wu, X., Wang, B., Xia, J.: Explicit symplectic multidimensional exponential
	fitting modified {R}unge-{K}utta-{N}ystr{\"o}m methods.
	\newblock BIT \textbf{52}(3), 773--795 (2012).
	\newblock \doi{10.1007/s10543-012-0379-z}
	
	\bibitem{Xu2022}
	Xu, Z., Cai, W., Hu, D., Wang, Y.: Exponential integrator preserving mass
	boundedness and energy conservation for nonlinear {S}chr\"odinger equation.
	\newblock Appl. Numer. Math. \textbf{173}, 308--328 (2022).
	\newblock \doi{10.1016/j.apnum.2021.12.007}
	
	\bibitem{ZK1974}
	Zakharov, V., Kuznetsov, E.: Three-dimensional solitons.
	\newblock Soviet Physics JETP \textbf{29}, 594--597 (1974)
	
\end{thebibliography}


{\color{black}\section*{Appendix}
	\setcounter{lemma}{1}
	\begin{lemma}
		For the semilinear system \eqref{semilinear HODE}, any exponential integrator of the form 
		\begin{align*}
			x_{n+1}=e^{hA} x_{n} +h \phi(hA)Q\hat{\nabla}U(x_{n},x_{n+1}),
		\end{align*}
		satisfies the following  equation
		\begin{align*}
			\frac{1}{2}x_{n+1}^\mathsf{T} Mx_{n+1}-\frac{1}{2}x_{n}^\mathsf{T} Mx_{n}+(x_{n+1}-x_{n})^T\hat{\nabla}U(x_{n},x_{n+1})=0,
		\end{align*}
		where  $\hat{\nabla}U(x_{n},x_{n+1})$ is a discretization of the true gradient $\nabla U(x_{n})$.
	\end{lemma}
	\begin{proof}
		We first assume that the matrix $M$ is not singular. Denoting by $V=hA$ and $\breve{\nabla}U=M^{-1}\hat{\nabla}U$ and using $\phi(V)V=e^{V}-I$, we get the following equation after replacing $x_{n+1}$ by $e^{hA} x_{n} +h \phi(hA)Q\hat{\nabla}U(x_{n},x_{n+1})$
		\begin{align*}
			&\frac{1}{2} x_{n+1}^\mathsf{T} M x_{n+1} -\frac{1}{2} {x_{n}}^{\mathsf{T}} M x_{n}+(x_{n+1} - x_n)^\mathsf{T} \hat{\nabla}U(x_n, x_{n+1})\\
			&=\frac{1}{2} {x_{n}}^{\mathsf{T}}({e^{V}}^{\mathsf{T}} M e^{V}-M)x_{n}+\frac{1}{2}h^2 \hat{\nabla}U^\mathsf{T} Q^\mathsf{T} \phi(V)^\mathsf{T} M \phi(V) Q\hat{\nabla}U\\
			&+ h {x_{n}}^{\mathsf{T}} {e^{V}}^{\mathsf{T}} M \phi(V) Q \hat{\nabla} U+{x_{n}}^\mathsf{T} ({e^{V}}^{\mathsf{T}}- I) \hat{\nabla}U + h{\hat{\nabla}U}^\mathsf{T} Q^\mathsf{T} \phi(V)^\mathsf{T}  \hat{\nabla}U\\
			&= \frac{1}{2} {x_{n}}^{\mathsf{T}} ({e^{V}}^{\mathsf{T}} M e^{V}-M)x_{n} 
			+\frac{1}{2}\breve{\nabla}U^\mathsf{T} (e^{V} - I)^\mathsf{T} M (e^{V} - I) \breve{\nabla}U\\
			&+{x_{n}}^\mathsf{T} ({e^{V}}^{\mathsf{T}} M - M)\breve{\nabla}U + {\breve{\nabla}U}^\mathsf{T} ({e^{V}}^{\mathsf{T}} M - M) \breve{\nabla}U+ {x_{n}}^{\mathsf{T}}{e^{V}}^{\mathsf{T}} M (e^{V} - I) \breve{\nabla}U\\
			&=\frac{1}{2}(x_{n}+\breve{\nabla}U)^\mathsf{T}({e^{V}}^{\mathsf{T}} M e^{V}-M)(x_{n}+\breve{\nabla}U)+\frac{1}{2}\breve{\nabla}U^\mathsf{T}({e^{V}}^{\mathsf{T}} M - Me^{V})\breve{\nabla}U\\
			&=0.
		\end{align*}
		The last equation follows from the fact that ${e^{V}}^{\mathsf{T}} M e^{V}-M$ and ${e^{V}}^{\mathsf{T}} M - M e^{V}$ are skew symmetric. When $M$ is singular, we follow the idea in \cite{Wuxinyuan2016}, introducing a series of symmetric and nonsingular matrices $M_\delta$ which converges to $M$ when $\delta$ approaches zero to complete the proof.
	\end{proof}
	
	\setcounter{corollary}{0}
	\begin{corollary}
		Let the potential $U$ be a general cubic polynomial with homogeneous decomposition as
		\[
		U(x)=U_3(x)+U_2(x)+U_1(x)+U_0,
		\]
		where $U_j$ is homogeneous of degree $j$.
		Then the EKahan scheme \eqref{EKahan scheme} satisfies
		\[
		H_{n+1}-H_n = U_3(x_{n+1}-x_n),
		\]
		showing that the stepwise energy defect depends only on the cubic component.
	\end{corollary}
	\begin{proof}
		By Lemma~\ref{lemma-energy} with $\hat{\nabla}U=\nabla_{\mathrm K}U$ and linearity of $\nabla_{\mathrm K}$ in $U$ and the linearity of $\nabla_{\mathrm K}$ regarding $U$, it suffices to treat each homogeneous part.
		For the homogeneous cubic part $U_3$, Theorem~\ref{theorem-energy} yields
		\[
		U_3(x_{n+1})-U_3(x_n)-(x_{n+1}-x_n)^\top \nabla_{\mathrm K}U_3(x_n,x_{n+1})
		=U_3(x_{n+1}-x_n).
		\]
		For $U_2$ and $U_1$, the Kahan polarization coincides with the exact discrete gradient identity
		\[
		(x_{n+1}-x_n)^\top \nabla_{\mathrm K}U_j(x_n,x_{n+1})=U_j(x_{n+1})-U_j(x_n),
		\qquad j=1,2,
		\]
		and $U_0$ contributes nothing.Therefore, only the cubic component contributes to the energy drift.
	\end{proof}

	\begin{lemma}
		Let $U$ be a polynomial of degree $k+2$. Then for each $n$ there exist a matrix $B_n\in\mathbb{R}^{d\times d}$ and a vector $\ell_n\in\mathbb{R}^d$, depending only on $(x_n,\dots,x_{n+k-1})$, such that
		\begin{equation}\label{eq:affine-polar appendix}
			\nabla_{\mathrm{K}}U(x_n,\dots,x_{n+k}) = B_n x_{n+k} + \ell_n.
		\end{equation}
		Consequently, \eqref{EKahan scheme-high-order} is equivalent to the linear system
		\begin{equation}\label{eq:linear-system-highorder appendix}
			\Bigl(I-kh\,\phi(khQM)\,Q B_n\Bigr)x_{n+k}
			= e^{khQM}x_n + kh\,\phi(khQM)\,Q\ell_n .
		\end{equation}
		In particular, if $U$ is homogeneous, then $\ell_n=0$.
	\end{lemma}
	\begin{proof}
		Fix $(x_n,\dots,x_{n+k-1})$ and regard $\nabla_{\mathrm{K}}U(x_n,\dots,x_{n+k})$ in~\eqref{function-polarization}
		as a function of the last argument $x_{n+k}$.
		Since $U$ is a polynomial, $\nabla U$ is a polynomial vector field, and the polarization formula~\eqref{function-polarization}
		is obtained from the symmetric polarization of each homogeneous component of $\nabla U$.
		Therefore $\nabla_{\mathrm{K}}U(x_n,\dots,x_{n+k})$ is \emph{linear} in $x_{n+k}$, and we may write
		\[
		\ell_n:=\nabla_{\mathrm{K}}U(x_n,\dots,x_{n+k-1},0),\qquad
		B_n x := \nabla_{\mathrm{K}}U(x_n,\dots,x_{n+k-1},x)-\ell_n,
		\]
		which yields~\eqref{eq:affine-polar appendix}. Substituting~\eqref{eq:affine-polar appendix} into~\eqref{EKahan scheme-high-order}
		gives~\eqref{eq:linear-system-highorder appendix}.
		Finally, if $U$ is homogeneous of degree $k+2$, then $\nabla U$ is homogeneous of degree $k+1$ and its polarization is $(k\!+\!1)$-linear,
		so $\nabla_{\mathrm{K}}U(x_n,\dots,x_{n+k-1},0)=0$, i.e., $\ell_n=0$.
	\end{proof}

	\begin{lemma}
		Let $A\in\mathbb R^{d\times d}$ and $g\in C^2([0,kh];\mathbb R^d)$. For all $t\in[0,kh]$, we claim
		\[
		\int_0^{kh} e^{(kh-s)A}g(s)\,ds
		=kh\,\phi(khA)\,g(kh/2)+R,
		\qquad \|R\|\le C h^3,
		\]
		where $C$ depends on $k$, $M_E$, $\|A\|$, and $\max\limits_{0\le s\le kh}\bigl(\|g(s)\|+\|g'(s)\|+\|g''(s)\|\bigr)$.
	\end{lemma}
	
	\begin{proof}
		Set $s_*=kh/2$ and $F(s)=e^{(kh-s)A}g(s)$. Then $F\in C^2$ and by the classical midpoint rule with Peano remainder, we have
		\begin{equation}\label{eq:midpoint-remainder}
			\int_0^{kh} F(s)\,ds = kh\,F(s_*) + R_1,
			\qquad
			\|R_1\|\le \frac{(kh)^3}{24}\max_{0\le s\le kh}\|F''(s)\|.
		\end{equation}
		A direct differentiation gives
		\[
		F''(s)=e^{(kh-s)A}\Bigl(A^2 g(s)-2A g'(s)+g''(s)\Bigr),
		\]
		hence
		\[
		\max_{0\le s\le kh}\|F''(s)\|
		\le M_E\Bigl(\|A\|^2\max\|g\|+2\|A\|\max\|g'\|+\max\|g''\|\Bigr).
		\]
		Therefore, there exist $C_1$ such that $\|R_1\|\le C_1 h^3$.
		
		Defining a function $D(z)=e^{z/2}-\phi(z)$, we have  $D(z)=-\frac{z^2}{24}+O(z^3)$ as $z\to 0$.
		Thus, 
		there exists $C_2$, depending on $k$, $M_E$, $\|A\|$, such that
		\[
		\|D(khA)\|\le C_2 h^2;
		\]
		and therefore, there exists $C_3$ such that
		\begin{equation*}
			kh\,e^{(kh/2)A}g(s_*) = kh\,\phi(khA)g(s_*) + kh\,D(khA)g(s_*)
			= kh\,\phi(khA)g(kh/2) + R_2,
		\end{equation*}
		and $ \|R_2\|\le C_3 h^3$. Combining with \eqref{eq:midpoint-remainder} yields the claim with $R=R_1+R_2$, and $C=C_1+C_3$.
	\end{proof}
	
	\begin{lemma}
		Assume $U$ is homogeneous of degree $k+2$ and let $x(t)$ be a smooth solution on $[t_n,t_n+kh]$.
		Set $\tilde x_{n+i}:=x(t_n+ih)$ and $t_{n+\frac{k}{2}}:=t_n+\frac{k}{2}h$.
		Then there exists $C>0$ (depending on $U$ and bounds of $x(t)$ on the macro-interval) such that
		\[
		\bigl\|\nabla_{\mathrm K}U(\tilde x_n,\ldots,\tilde x_{n+k})
		-\nabla U\bigl(x(t_{n+\frac{k}{2}})\bigr)\bigr\|
		\le C h^2 .
		\]
	\end{lemma}
	\begin{proof}
		Let $\bar x:=x(t_{n+\frac{k}{2}})$ and write $\tilde x_{n+i}=\bar x+\delta_i$.
		A Taylor expansion of $x(t)$ around $t_{n+\frac{k}{2}}$ gives
		\begin{equation*}
			\delta_i=(i-\tfrac{k}{2})h\,\dot x(t_{n+\frac{k}{2}})+ O(h^2),
			\qquad i=0,\ldots,k,
		\end{equation*}
		hence $\|\delta_i\|= O(h)$ and moreover
		\[
		\sum_{i=0}^k \delta_i= O(h^2),
		\quad\text{since}\quad \sum_{i=0}^k (i-\tfrac{k}{2})=0.
		\]
		For homogeneous $U$, the polarized gradient is symmetric and multilinear in its $(k+1)$ arguments; in particular,
		$\nabla_{\mathrm K}U(\bar x,\ldots,\bar x)=\nabla U(\bar x)$.
		Using multilinearity around the diagonal point $(\bar x,\ldots,\bar x)$, we expand
		\[
		\nabla_{\mathrm K}U(\bar x+\delta_0,\ldots,\bar x+\delta_k)
		=\nabla_{\mathrm K}U(\bar x,\ldots,\bar x)
		+\sum_{i=0}^k \mathcal L(\delta_i)
		+ O(h^2),
		\]
		where $\mathcal L(\cdot)$ denotes the same linear map (by symmetry) obtained by replacing one argument $\bar x$
		with the increment. Consequently,
		\[
		\Bigl\|\sum_{i=0}^k \mathcal L(\delta_i)\Bigr\|
		\le C_1\Bigl\|\sum_{i=0}^k \delta_i\Bigr\|
		= O(h^2),
		\]
		and the remaining terms contain at least two increments and are therefore also $ O(h^2)$.
		Combining these estimates yields
		\[
		\nabla_{\mathrm K}U(\tilde x_n,\ldots,\tilde x_{n+k})
		=\nabla U(\bar x)+ O(h^2),
		\]
		which concludes the proof.
	\end{proof}
	
	\begin{lemma}
		Assume the exact solution stays in a bounded set $\mathcal D$ on $[0,T]$ and $U$ is polynomial.
		Then there exist $h_0>0$ and $C>0$ such that for all $0<h\le h_0$ the $k$-step EKahan update
		defines a unique one-step map $Y_{n+1}=\Phi_h(Y_n)$ on $\mathcal D^k$, and
		\[
		\|\Phi_h(Y)-\Phi_h(Z)\|_\infty\le (1+Ch)\,\|Y-Z\|_\infty,
		\qquad \forall\,Y,Z\in\mathcal D^k,
		\]
		where $\|Y\|_\infty:=\max\limits_{0\le j\le k-1}\|y_j\|$ for $Y=(y_0,\ldots,y_{k-1})\in(\mathbb R^d)^k$.
	\end{lemma}
	\begin{proof}
		From Theorem \ref{thm-highorder-wellposed}, the update from $Y_n$ to $Y_{n+1}$ can be uniquely defined via solving \eqref{eq:linear-system-recall}.  For given $Y$ and $Z$ in $\mathcal D^k$, we claim ,  we claim
		\[
		\|\Phi_h(Y)-\Phi_h(Z)\|_\infty
		=\max\Bigl\{\max_{1\le j\le k-1}\|y_j-z_j\|,\ \|y_k-z_k\|\Bigr\}
		\le (1+\widetilde C h)\,\|Y-Z\|_\infty, 
		\]
		based on the estimation in \eqref{high-esti-one-vec}.
	\end{proof}
	\begin{corollary}
		For the semilinear system \eqref{semilinear HODE}, any exponential integrator of the form
		\[
		x_{n+k}=e^{khA} x_{n} +kh \phi(khA)\,Q\,\hat{\nabla}U(x_{n},x_{n+1},\cdots,x_{n+k}),
		\]
		where  $\hat{\nabla}U(x_{n},x_{n+1},\cdots,x_{n+k})$ is a discretization of the true gradient $\nabla U(x_n)$, satisfies
		\[
		\frac{1}{2k}x_{n+k}^\mathsf{T} Mx_{n+k}-\frac{1}{2k}x_{n}^\mathsf{T} Mx_{n}
		+\frac{1}{k}(x_{n+k}-x_{n})^\mathsf{T} \hat{\nabla}U(x_{n},x_{n+1},\cdots,x_{n+k})=0.
		\]
	\end{corollary}
	\begin{proof}
		Assume first that $M$ is nonsingular. Set $V=khA$ and define $\breve\nabla U=M^{-1}\hat\nabla U$.
		Using $\phi(V)V=(e^{V}-I)$ and $A=QM$, we obtain
		\[
		kh\,\phi(V)\,Q\,\hat\nabla U
		=(e^{V}-I)\breve\nabla U,
		\qquad
		x_{n+k}-x_n=(e^{V}-I)(x_n+\breve\nabla U).
		\]
		A direct expansion gives
		\begin{align*}
			&\frac{1}{2k}x_{n+k}^TMx_{n+k}-\frac{1}{2k}x_n^TMx_n+\frac{1}{k}(x_{n+k}-x_n)^T\hat\nabla U \\
			&=\frac{1}{2k}(x_n+\breve\nabla U)^T\bigl(e^{V^T}Me^V-M\bigr)(x_n+\breve\nabla U)
			+\frac{1}{2k}\breve\nabla U^T\bigl(e^{V^T}M-Me^V\bigr)\breve\nabla U,
		\end{align*}
		which vanishes since $V^\mathsf{T}M+MV=0$ implies $e^{V^\mathsf{T}}Me^V=M$ and $e^{V^T}M-Me^V$ is skew-symmetric. The singular case follows the same as in Theorem \ref{theorem-energy}.
	\end{proof}
	
	\begin{corollary}
		For a homogeneous function $U(x)$ of degree $k+2$, EKahan scheme \eqref{EKahan scheme-high-order} satisfies the stepwise identity
		\begin{align*}
			H_{n+1}-H_{n}=&\bar{U}(x_{n+1},x_{n+2},\dots,x_{n+k-1},x_{n+k},x_{n+k},x_{n+1}-x_{n})\\
			&+\bar{U}(x_{n},x_{n+1},\dots,x_{n+k-2},x_{n+k-1},x_{n},x_{n+k}-x_{n+k-1})\\
			&-2\bar{U}\Bigl(x_{n},x_{n+1},\dots,x_{n+k-1},x_{n+k},\frac{x_{n+k}-x_{n}}{k}\Bigr),
		\end{align*}
		where
		\[
		H_{n}=\frac{1}{2k}\sum_{i=0}^{k-1} x_{n+i}^\mathsf{T} Mx_{n+i}
		+\bar{U}(x_{n},x_{n+1},\dots,x_{n+k-2},x_{n+k-1},x_{n},x_{n+k-1})
		\]
		is the discrete energy at $t=t_n$, and $\bar{U}(\cdot,\dots,\cdot)$ is a symmetric $(k+2)$-linear form associated with $U$ satisfying $\bar{U}(x,\dots,x)=U(x)$.
	\end{corollary}
	
	\begin{proof}
		Applying Corollary~\ref{corollary-high-order-energy} with $\hat\nabla U=\nabla_{\mathrm K}U(x_n,\ldots,x_{n+k})$ gives
		\begin{equation}\label{eq:quad-diff-high}
			\frac{1}{2k}x_{n+k}^TMx_{n+k}-\frac{1}{2k}x_n^TMx_n
			=-\frac{1}{k}(x_{n+k}-x_n)^T\nabla_{\mathrm K}U(x_n,\ldots,x_{n+k}).
		\end{equation}
		Recall that $U$ is homogeneous of degree $k+2$ and can be represented by a symmetric $(k+2)$-linear form $\bar U$
		such that $U(x)=\bar U(x,\ldots,x)$.
		Moreover, the homogeneity implies the directional derivative identity
		\begin{equation*}
			u^T\nabla U(z)=(k+2)\,\bar U(u,z,\ldots,z)
			\qquad \text{for all } u,z\in\mathbb R^d.
		\end{equation*}
		Since $\nabla_{\mathrm K}U$ is the $(k+1)$-linear polarization of $\nabla U$ and polarization commutes with linear
		functionals, it follows that for all $u\in\mathbb R^d$,
		\begin{equation*}
			u^T\nabla_{\mathrm K}U(x_n,\ldots,x_{n+k})=(k+2)\,\bar U(u,x_n,\ldots,x_{n+k}).
		\end{equation*}
		Taking $u:=x_{n+k}-x_n$ and using multilinearity of $\bar U$ yields
		\begin{align}\label{eq:key-polar-identity}
			\frac{1}{k}u^T\nabla_{\mathrm K}U(x_n,\ldots,x_{n+k})
			=2\bar U\!\Bigl(x_n,x_{n+1},\ldots,x_{n+k},\frac{u}{k}\Bigr)
			+\bar U(x_n,x_{n+1},\ldots,x_{n+k},u).
		\end{align}
		By the definition of the discrete energy $H_n$ in the statement,
		\begin{align*}
			H_{n+1}-H_n
			&=\frac{1}{2k}(x_{n+k}^TMx_{n+k}-x_n^TMx_n)
			+\bar U(x_{n+1},\ldots,x_{n+k},x_{n+1},x_{n+k})\\
			&-\bar U(x_n,\ldots,x_{n+k-1},x_n,x_{n+k-1}).
		\end{align*}
		Combining \eqref{eq:quad-diff-high} and \eqref{eq:key-polar-identity} and using the multilinearity
		and symmetry of $\bar U$ yields the stated increment formula.
\end{proof}}

\end{document}